\documentclass[9pt]{amsart}
\usepackage[square,compress,comma, numbers,sort]{natbib}
\usepackage[colorlinks=true, citecolor=red, linkcolor=blue]{hyperref}
\usepackage{amsfonts,mathtools}
\usepackage{xparse}
\usepackage{amssymb} 
\usepackage[titletoc,title]{appendix}
\usepackage{color}
\usepackage{mathrsfs}   
\usepackage[shortlabels]{enumitem}
\usepackage{cleveref}

\def\drz{\lambda_\alpha(dz)}
\def\drzz{\alpha z^{\alpha- 1} dz}

\def\drz{\alpha z^{-\alpha- 1} dz}

\def\drzz{\alpha z^{\alpha- 1} dz}

\def\wH{ { \mathfrak{H}}}
\def\wwH{{\wH_\star}}

\def\Pspace{ \parens{\Omega,\mathscr{F},\mathbb{P}}}

\newcommand\Dset{\mathrm{D}(\TT, \R^d)}
\newcommand\DDset{D_c}

\DeclarePairedDelimiter{\parens}()
\def\Pspace{ \parens{\Omega,\mathscr{F},\mathbb{P}}}

 \def\borel#1{\mathscr{B}(#1)}

\def\ve{\varepsilon}

\def\Ife{\mathcal{I}_{\EI,fe}}  
\def\Ifasup{ \mathcal{I}_{\EI,argsup}}   
\def\II{\mathfrak{I}}
\def\bobo{\mathfrak{\vartheta}}

\def\SZq{ \mathcal{S}_{\EIE} (\vk Z)}
\def\STq{ \mathcal{S}_{\EIE} (\vk \Theta)}

\def\JJqx{ \mathcal{B}_{\EIE,\tau} (\vk Y )}

\allowdisplaybreaks[4]

\newcommand{\abs}[1]{\left\lvert #1 \right\rvert}

\DeclarePairedDelimiterXPP\pk[1]{\mathbb{P}}\{ \}{}{ #1}
\DeclarePairedDelimiterXPP\E[1]{\mathbb{E}}\{ \}{}{	#1}

\DeclarePairedDelimiterXPP\ind[1]{\mathbb{I}}( ){}{	#1}

\def\FRE{\mbox{Fr\'{e}chet }}

\def\toas{\overset{a.s.}\rightarrow}
\def\toprob{\overset{p}\rightarrow}

\def\wHO{ \wH_{0}}
\def\wHA{ \wH_{\alpha}}

\newcommand{\norm}[1]{\lVert  #1 \rVert }

\definecolor{c20}{rgb}{0.,0.7,0.}
\definecolor{c30}{rgb}{0.,0.,1.}
\definecolor{c40}{rgb}{1,0.1,0.7}
\definecolor{c50}{rgb}{1,0,0}
\definecolor{c60}{rgb}{1,0.9,0.1}
\definecolor{c70}{rgb}{0.50,1.00,0.00}

\def\eB#1{#1}
\def\cEE#1{#1} 

\def\cEp#1{#1}

\def\cER#1{#1}

\def\cED#1{#1}

\numberwithin{equation}{section}

\newtheorem{theo}{Theorem}[section]
\newtheorem{sat}[theo]{Proposition}
\newtheorem{de}[theo]{Definition}
\newtheorem{lem}[theo]{Lemma}

\newtheorem{example}[theo]{Example}
\newtheorem{korr}[theo]{Corollary}
\newtheorem{remark}[theo]{Remark}

\numberwithin{equation}{section}

\newtheorem{condition}{Condition}[section]

\newcommand{\prooftheo}[1]{ \textsc{Proof of Theorem} \ref{#1} }
\newcommand{\proofprop}[1]{\textsc{Proof of Proposition} \ref{#1}}
\newcommand{\prooflem}[1]{\textsc{Proof of Lemma} \ref{#1}}

\newcommand{\QED}{\hfill $\Box$}
\newcommand{\COM}[1]{}
\def\IF{\infty}
\newcommand{\R}{\mathbb{R}}
\newcommand{\inr}{\in \R}

\newcommand{\BQN}{\begin{eqnarray}}
\newcommand{\EQN}{\end{eqnarray}}
\newcommand{\BQNY}{\begin{eqnarray*}}
\newcommand{\EQNY}{\end{eqnarray*}}
\def\ldot{, \ldots,}

\newcommand{\limit}[1]{\lim_{#1 \to   \infty}}

\def\todis{\overset{d}\rightarrow}

\newcommand{\vk}[1]{#1}

\def\cadlag{c\`adl\`ag}

\def\ifft{iff}
\def\shiftinvariant{shift-invariant}
\def\shiftinvariance{shift-invariance }
\def\Shiftinvariance{Shift-invariance }

\def\ttt{\\&=&}

\def\bqny#1{ \begin{eqnarray*} #1 \end{eqnarray*}}
\def\bqn#1{ \begin{eqnarray} #1 \end{eqnarray}}

\newcommand{\BS}{\begin{sat}}
\newcommand{\ES}{\end{sat}}
\newcommand{\BT}{\begin{theo}}
\newcommand{\ET}{\end{theo}}
\newcommand{\BK}{\begin{korr}}
\newcommand{\EK}{\end{korr}}

\newcommand{\BEX}{\begin{example}}
\newcommand{\EEX}{\end{example}}

\newcommand{\BD}{\begin{de}}
\newcommand{\ED}{\end{de}}

\newcommand{\BRM}{\begin{remark}}
\newcommand{\ERM}{\end{remark}}

\newcommand{\BEL}{\begin{lem}}
\newcommand{\EEL}{\end{lem}}
\def\mj{\mathcal{J}}

\def\JJ{{\mathcal{B}_\tau}(\vk Y)}

\def\Z{\mathbb{Z}}
\def\inn{\in \mathbb{N}}

\def\MM{\mathcal{M}(\vk \Theta)}

\def \SZ{ \mathcal{S}( {\vk Z})}

\def\ST{ \mathcal{S}( {\vk \Theta})}

\def\Pspace{ \parens{\Omega,\mathscr{F},\mathbb{P}}}
\def\TTd{\EIE }
\def\kk{\mathcal{ K}_{\alpha}[\vk Z]}

\def\TTT{\mathbb{T} }
\def\TT{\mathcal{T} }
\def\intT{\int_{\TT}}

\def\ZT{\widetilde{ \vk Z}}   

\def\kK{\mathscr{C}_{\kappa} [\vk Z]}

\def\EIE{ \mathcal{L}}
\def\EI{ \mathcal{L}}

\def\STYq{\mathcal{S}_\EI(\vk Y)}

\def\kKY{\mathbb{K}_\alpha[\vk Z]}
\def\WZ{ \widetilde{\vk Z} }

\def\TTT{\TT}

\def\diad{{\mathbb{T}_0}}

\def\intD#1{\sum_{#1\in \diad}}

\def\kk{\mathcal{C}_{\kappa} [\vk Z ]}

\def\pclas{\mathcal{P}_\diad}
\def\intDD{\int_{\diad}}
\def\pclasT{\mathcal{P}_\TT}

\def\AA{\mathscr{D}}

\def\Hh{ \mathfrak{H}}

\def\clasU{\mathbb{U}}

\DeclarePairedDelimiterXPP\ETT[1]{\widetilde{\mathbb{E}}}\{ \}{}{	#1}

\def\SY{\mathcal{S}(\vk Y)}

\def\vkT{\cER{\vk \Theta_\kappa}}
\def\tX{\cER{\kappa(\vk Z)}}
\def\vkT{\cER{\vk \Theta}}
\def\SCD{\mathfrak{W}_{\kappa}}

\begin{document}

	\title{Shift-invariant homogeneous classes of random fields}

	\author{Enkelejd  Hashorva}
	\address{Enkelejd Hashorva, Department of Actuarial Science,  
		University of Lausanne,\\
		UNIL-Chamberonne, 1015 Lausanne, Switzerland}
\email{Enkelejd.Hashorva@unil.ch}

	\bigskip
	
	\date{\today}
	\maketitle

 \begin{quote} 
 {\bf Abstract:} 
Given an $\R^d$-valued  random field (rf) $\vk Z(t),t\in \TT$ and an $\alpha$-homogeneous  mapping $\kappa$  we  define  the corresponding equivalent class of rf's (denoted by $\kK$) which include representers of the same tail measure $\nu_Z$. When $\TT$ is an additive group, tractable equivalent classes of interest are the shift-invariant ones, which contain in particular all  independent random shifts of $\vk Z$.  This contribution is mainly  concerned with   the investigation  of the probabilistic properties of shift-invariant $\kK$'s. Important objects introduced in our setting are tail and spectral tail rf's.  Further, the class of universal  maps $\clasU$ acting on elements of $\kK$ turns out to be crucial for properties of functionals of $\vk Z$.  Applications of our findings concern   max-stable and symmetric $\alpha$-stable  rf's,  their maximal indices as well as their random shift-representations.
 \end{quote}


\section{Introduction}
Homogeneous  functionals (hf's) play a crucial role in the study of max-stable and symmetric $\alpha$-stable random fields (rf's) as well as tail measures and tail rf's. This fact has been  exploited in  \cite{MolchanovBE}, see also the later contributions \cite{Htilt,Stb,kremer,Hrovje, kulik:soulier:2020,klem,BP,PH2020,KumeE}. It is simpler  to demonstrate the importance of hf's in the context of max-stable rf's as we show next. 
Fix below  $d,l$ two positive integers, $\alpha>0$ and a norm 
  $\norm{\cdot}$ on $\R^d$.  Let $\vk X(t),t\in \TT $  be an $\R^d$-valued max-stable rf with $\alpha$-\FRE max-stable marginal distributions having almost surely (a.s.) \cadlag\ sample paths being   defined on a complete probability space $\Pspace$. In this contribution we consider $\TT=\R^l$ or \cER{$\TT\subset \Z^l$} assumed to be  non-empty. In the light of  de Haan characterisation,  see e.g., \cite{deHaan,stoev2005extremal, MolchanovBE,MolchanovSPA}  we have the following representation (in distribution)
\bqn{\label{eq1}
	\vk X(t) =  \max_{i\ge 1} \Gamma_i^{-1/\alpha} \vk Z^{(i)}(t), \quad t\in \TT,  
}
where $\Gamma_i= \sum_{k=1}^i \mathcal{E}_k$ with $\mathcal{E}_k, k\ge 1$ unit iid exponential random variables (rv's) independent of $\vk Z^{(i)}$'s, which are independent copies of a $d$-dimensional rf $\vk Z(t)=(Z_1(t) \ldot Z_d(t)),t\in \TT$ with non-negative components having a.s.\ \cadlag\ sample paths such that $\E{ \norm{\vk Z(t)}^\alpha}< \IF$ for all $t\in \TT$. 
 The rf $\vk Z$ is referred to as the representer  of $\vk X$.  It is well-known that for all $t_i \in \TT, 
\vk x_i\in (0,\IF)^d, i=1 \ldot n, \gamma>0$   
\bqn{\label{kode}
	(\pk{ \vk X(t_i) \le \vk x_i, i=1 \ldot  n})^\gamma&=& 
	\pk{ \vk X(t_i) \le \vk x_i/\gamma ^{1/\alpha}, i=1 \ldot  n} \notag\\
	& =& \exp \Biggl( - \gamma \E*{ \max_{ 1 \le i\le d, 1 \le j \le n   }  Z_i^\alpha(t_j)/x_{ij}^\alpha    }  \Biggr).
}

We shall consider next rf's $\vk Z(t),\ZT(t),t\in \TT$ that satisfy (the backshift  opreration is $B^{-t} f= f(t)$) 
 \bqn{
 	\label{ZY} 
 	\pk*{ \sup_{t\in \TTT }\kappa(B^{-t} \vk Z) >0}= 	\pk*{ \sup_{t\in \TTT } \kappa(B^{-t} \ZT)>0}=1,
 	\quad \kappa(B^{-t} f)= \norm{f(t)}^\alpha,
 }  
where $\Dset$ denotes the space of functions $f: \TT\mapsto \R^d$  equipped with  the product (cylindrical) $\sigma$-field   $\AA$ and    
$$\sup_{t\in A} \kappa (B^t f):= \max_{t\in A \cap \diad} \kappa(B^t f), \quad A \subset \TT, f\in \Dset. 
$$
 Throughout this paper $\diad$ consists of all $t\in \R^l\cap \TT$ which have rational coordinates.\\
  In view of \cite[Lem 7.1]{HBernulli} the requirement \eqref{ZY} is  no restriction for a given representer $Z$ of $\vk X$. 
 Suppose hereafter without loss of generality that 
\bqn{ \label{lBo}
	\E*{\kappa(\vk{Z})  }\in (0,\IF).
}
As shown in \cite[Prop 2.1]{KumeE} given two  representers $\vk Z$ and $\ZT$ of the max-stable rf $\vk X$
\bqn{ \label{boll} 
	\E{\kappa(B^{-h} \ZT) F( \ZT /\kappa(B^{-h} \ZT)^{1/\alpha} )} = \E{\kappa(B^{-h} \vk Z) F( \vk Z /\kappa(B^{-h} \cED{\vk Z})^{1/\alpha} )}, 
	\quad  \forall F\in \mathcal{H}, \forall h\in \TT 
}
holds 
with $\mathcal{H}$ the class  of all  $\AA/\borel{(-\IF,\IF]}$-measurable maps $F: \Dset  \mapsto E$ where $E=\R$ or $E=[0,\IF]$.  Here   $\borel{V}$ stands for the Borel $\sigma$-field of a topological space $V$.\\
Moreover, also the converse is valid, i.e., any non-negative $\R^d$-valued rf $\ZT$ with \cadlag\ sample paths   that satisfies \eqref{ZY} and  \eqref{boll} is a representer for $\vk X$. \\

\BD  If $\TT=\R^l$, then  $\DDset\subset \Dset$ consists of only \cadlag\ functions,  while for $\TT \subset \Z^l$ we set $\DDset=\Dset$. In both cases $\DDset$ is equipped with a metric $d_{\DDset}$ that turns it into a Polish space.
\label{def:dc}
\ED
The functional identity \eqref{boll} is a consequence of Balkemma's Lemma, see \cite[Lem 4.1]{HaanPickands} and \cite[Rem 3.3]{KumeE}. Moreover, it can be related to the regular variation of rf's. Specifically, in view of \cite{kulik:soulier:2020,PH2020,MartinE} where \cadlag\  rf's are considered, $\vk X$ is regularly varying with respect to the boundedness $\mathcal{B}_0$,  some positive sequence $a_n>0,n\inn$ and the  tail measure $\nu_{\vk Z}$ defined  on $\AA$ 
by 
\bqn{ \label{nuZ} 
	\nu_{\vk Z}[F]= \int_0^\IF \E{F(z \vk Z)} \alpha z^{-\alpha- 1} dz, \quad \forall F\in \mathcal{H}.
}	
Tail measures are first introduced in \cite[p.\ 159]{MR3561100} and subsequently studied in \cite{wao,klem,Hrovje,kulik:soulier:2020,PH2020,MartinE}.\\
 Applying  \Cref{lemZJ} in Appendix we obtain 
\bqny{
\limit{n} n \pk{ \kappa(B^h  \vk X) F(\vk X)> b_n}= 
\E{\kappa( B^h \vk Z) F(\vk Z)}< \IF, \quad\forall h\in \TT
}
for all 
$F: \DDset \mapsto [0,\IF) $ continuous, $F(0)=0$, and $0$-homogeneous with an appropriate choice of $b_n$'s. 
Since for another representer $\ZT$ of $\vk X$ necessarily $\nu_{\ZT}=\nu_{\vk Z}$,  \cED{then} \eqref{boll} follows. From \eqref{kode} and the assumption on the sample paths of $\vk X$ we obtain 
\begin{equation}\label{conditionC1}
	\E*{ \sup_{t\in [-c,c]^l  \cap \TT} \kappa( B^{-t} \vk{Z}) } < \IF, \quad \forall \cER{c}\in (0,\IF).
\end{equation}
In the light of \eqref{boll} and the regular variation considerations above,  it becomes clear that hf's play a crucial role in the characterisation of the representers of $\vk X$. \\
{Two representers $\vk Z$ and $\ZT$ of $\vk X$ are called in \cite{MolchanovBE,KumeE} max-zonoid equivalent, while in \cite{Domb14} and papers that refer to it, those are simply called equivalent. 
	A flow representation of a representer of $\vk X$ is derived 
	in \cite[Thm 1]{MR2453345}, see also \cite[p.1219]{MolchanovBE} for $d=1$.  In fact, flow representations are commonly discussed  for symmetric $\alpha$-stable processes based on results of Rosi\'nski \cite{rosi1}, see e.g., \cite{MolchanovBE,MR3561100,Roy}. Note that two representers of $\vk X$ can be also related directly to each other, see  
	\cite[Eq.\ (3.1)]{MolchanovBE} for the case $\alpha=d=1$.
}	

Indeed, the imposed restriction that both $\vk Z, \ZT$ have  non-negative components is not essential.  In connection with symmetric $\alpha$-stable rf's this fact is known from \cite[Thm 1.1.]{StoevWangSPL} and \cite[Thm 2.2]{MolchanovBE}.\\
 We drop therefore that  specific assumption in the discussion below and consider $\nu_{\vk Z}$ on $\AA$ for  $ \vk Z$ that can have negative components. Also for such $ \vk Z$, in view of 
\cite[Prop 3.6, Rem 3.11,(ii)]{MartinE} the functional identity \eqref{boll} for 
$\ZT$ satisfying \eqref{ZY} is equivalent with 
\bqn{\nu_{\vk Z}=\nu_{\ZT}. 
}
Hereafter $\kappa: \Dset \mapsto [0,\IF]$ belongs to $\mathcal{H}$ being further $\alpha$-homogeneous   i.e., 
$$\kappa(cf) =c^\alpha \kappa(f),\quad  \forall c>0,\forall f\in \Dset.$$ 
\BD $\SCD$ stands for the class of $\R^d$-valued  rf's $\vk V(t), t \in \TT$ defined on  a complete probability spaces $(\widetilde \Omega, \widetilde {\mathscr{F}}, \widetilde{\mathbb{P}})$ such that $\kappa(B^{-t} \vk V), t\in \TT$ is stochastically continuous. Moreover, $\kappa(B^t\vk V),t\in \TT$ is separable with separant $\diad$ and jointly measurable.
\ED 
Constructing a new representer $\ZT$ for a given max-stable rf $\vk X$ or for a given tail measure $\nu_{\vk Z}$ is an important topic discussed in numerous papers. The former task is for instance crucial for simulations, see e.g., \cite{KabExt} and the references therein. In this article, we shall mainly focus on the functional equation \eqref{boll} discussing the properties of representers of $\nu_{\vk Z}$ from $\SCD$ without making any specific reference to regular variation. In particular 
\begin{enumerate}[i)] 
	\item instead of a norm $\norm{\cdot}$ on $\R^d$, below we shall consider an $\alpha$-homogeneous and $\AA/\borel{[0,\IF]}$-measurable map  $\kappa: 
		\Dset \mapsto [0,\IF] $;
		\item  we shall only assume that the representers $\ZT$ belong to $\SCD$ when $\TT= \R^l$. 
\end{enumerate}

The representer of $\nu_{\vk Z}$ that belong to $\SCD$ define a class of rf's. Note that below $\ZT$ might be defined in another probability space than $Z$.  For notational simplicity  we write $\mathbb{P}$ and $\mathbb{E}$ instead of $\widetilde{\mathbb{P}}$ and $\widetilde{\mathbb{E}}$, respectively. 
 \BD 
  $\kK$ with representer $\vk Z$ 
consists of $\vk Z,\ZT \in \SCD$ that satisfy  \eqref{ZY}-\eqref{boll} and \eqref{conditionC1}.
\ED 
 The Brown-Resnick max-stable rf's were introduced for particular instances initially in 
  \cite{bro1977, eddy1980distribution}  and discussed in greater generality later in \cite{kab2009}. 
   We introduce below    the corresponding  $\alpha$-homogeneous classes.  

\begin{example}[Brown-Resnick $\kK$]
	For $\vk V(t)=(V_1(t) \ldot V_d(t)), t\in \TT$ a centered $\R^d$-valued Gaussian rf with a.s.\ continuous sample paths 
define 
	$$\vk Z(t)= ( \xi_1 e^{  W_1(t)} \ldot 
	\xi_d e^{  W_d(t)} ), \quad  W_i(t)= V_i(t) - \alpha Var(V_i(t))/2, \quad 1 \le i \le d, t\in \TT, $$
	where  $\xi_1 \ldot \xi_d$ are rv's that take values $\pm 1$  being further independent of any other random element.   Taking  
	$$\norm{\vk x}_\alpha= (\sum_{i=1}^d \abs{x_i}^\alpha /d)^{1/\alpha},\quad \vk x\inr^d,$$
	 then  $ \E{ \norm{\vk Z(t)}^\alpha_\alpha }=1$ for all $t\in \TT$. In  light of 
	 \cite[Cor.\ 6.1]{MR3024389}
	$$d\E*{ \sup_{t\in [-c,c]^l} \norm{\vk Z(t)}^\alpha_\alpha}= \sum_{i=1}^d\E*{ \sup_{t\in [-c,c]^l} e^{ \alpha V_i(t) -  Var( \alpha V_i(t))/2 }} < \IF, 
	\ \ \forall c>0$$
	and thus     $\vk Z \in \SCD$ and it  satisfies    
	 \eqref{ZY}, \eqref{lBo}, \eqref{conditionC1} with $\kappa(f)= \norm{f(0)}^\alpha_\alpha$. Hence we can define the corresponding $\kK$. 
Write $\gamma$ for the pseudo-variogram matrix-valued function of $\vk V$ with  components  $$\gamma_{ij}(s,t)= Var(V_i(s)-V_j(t))/2,\quad s,t\in \TT, 1 \le i \le d, 1 \le j \le d.$$
In view of \Cref{examp:5B} \cED{below} $\gamma$ uniquely   determines this class.
	\label{examp:5}
\end{example}	 

The Brown-Resnick-L\'evy max-stable  processes have been introduced in \cite{StoevSPA}, see \cite{eng2014d,Htilt} for further results. The corresponding $\alpha$-homogeneous class of rf's is defined next.  
\begin{example}[Brown-Resnick-L\'evy $\kK$]\label{exa:lev} Let $W_{ij}(s),s\ge 0,1\le i\le l, 1\le j\le d$ be L\'evy processes with Laplace exponent $\psi_{ij}(\theta)=\ln \E{ e^{\theta W_{ij}(1)}}$ such that $\psi_{ij}(\alpha)=0$ for some $\alpha>0, 1 \le i\le l, 1 \le j\le d$. Write $W_{ij}^{(\alpha)}(s),s\ge 0$ for the exponentially tilted L\'evy process with Laplace exponent $\psi_{ij}(\alpha+\theta)$. Suppose further that $W_{ij}, W_{ij}^{(\alpha)}$'s are all independent with \cadlag\ sample paths defined in the same probability space and set 
	$$  Z_j(t)= \xi_j\prod_{i=1}^l e^{ \ind{t_i\ge 0}W_{ij}(t_i)  - \ind{t_i<0}  W_{ij}^{(\alpha)}((-t_i)-)}, \quad t=(t_1 \ldot t_l),$$
	with  $\xi_i$'s as in \Cref{examp:5}.  Taking $\kappa(f)= \norm{f(0)}^\alpha_\alpha$ it follows easily that $\kK$ is well-defined. 
\end{example}
In applications,  of particular interest are max-stable, symmetric $\alpha$-stable, and  tail measures that are 
shift-invariant. A unified approach in their study is the investigation of  shift-invariant  $\alpha$-homogeneous classes of rf's, which we define below. 
 Whenever we consider the shift-invariance, we shall assume that $\TT$ is an additive group and hence $B^h f, f\in \Dset$ is well-defined.

 \BD  $\kK$ is shift-invariant if for some $\ZT \in \kK$ we have 
\bqn{\label{defShift}
	B^h \ZT \in \kK, \quad \forall h\in \TT.
}	
\ED   
\cER{By the above  definition and  \eqref{boll} if $\kK$ is  shift-invariant, then  
\bqn{ \quad \quad  
	\E*{ \kappa(B^{-h}\vk Z^\star) F( \vk Z^\star )}&=&
	\E*{ \kappa( B^{-h}\ZT)  F( \ZT )}
	=\E{ \kappa( \WZ ) F(B^h \WZ) }\notag \\
	&=&\E*{ \kappa(  \vk Z^\star )  F( B^h\vk Z^\star )}, \label{tcfN0} 
\quad   \forall F\in \mathcal{H}_{0},\forall h\in \TT, \forall \vk Z^\star \in \kK 
} 
since when $F\in \mathcal{H}_{0}$ also $G(f)=F(B^h f), f\in \Dset$   belongs to $\mathcal{H}_{0}$  
and thus \eqref{defShift} holds for all $\vk Z^\star \in \kK$.} 
Clearly,  if $\vk Z$ is stationary, then $\kK$ is shift-invariant. More interesting  cases are shift-invariant $\kK$'s generated by some non-stationary $\vk Z$. Two prominent instances are  the Brown-Resnick $\kK$, see  \Cref{examp:5B}, while the \shiftinvariance of the Brown-Resnick-L\'evy $\kK$ follows easily from the definition and the \shiftinvariance of the corresponding max-stable rf, see \cite{eng2014d}. \\

It was shown in \cite{Htilt} for the case $d=1$ and $\vk Z$ being non-negative and in \cite{Hrovje,klem,kulik:soulier:2020,BP, PH2020,MartinE,Planic} for discrete or \cadlag\ $\vk Z$ that the functional equation \eqref{boll} is equivalent with the  time-change formula discovered in \cite{BojanS}. Therein it is  formulated in terms of the so called spectral tail rf $\vk \Theta$, see below \eqref{eqDo20}. As shown in  \cite{Hrovje} the time-change formula is equivalent to  \eqref{tYY}. Dropping the regular variation context of \cite{BojanS}, we introduce a more general object, again denoted by $\vk \Theta$, which as in \cite{klem} is labeled as the local rf.
\BD  \cER{Given a $\kK$   we shall define its   local rf $\vkT$  as the rf $\vk Z/\tX^{1/\alpha}$} under the probability law $\widehat{ \mathbb{P}}$, where (recall that we assume \eqref{lBo}) 
\begin{equation}\label{minist}
  \widehat {\mathbb{P}}(A)= \frac{1}{\E*{\tX  }}\E{   \tX \ind{ A}}, 
  \  \
\quad \forall A \in \mathscr{F},
\end{equation} 
with $\ind{A}$   the indicator function of the set $A$. 
\ED 
Brief organisation of the paper: \Cref{sec:prim} introduces lattices, some classes of maps and  our  notation. 
Extension of \eqref{boll} is  discussed  in  \Cref{sec:shift} followed by a section on the characterisations of shift-invariance $\kK$.  Both spectral tail and tail rf's are introduced  in \Cref{secTailSpectral}. 
Universal maps and their relations with   shift-involutions and anchoring maps 
are introduced in \Cref{sec:uni}. 
 Implications of our findings are briefly considered in \Cref{sec:impli}. All proofs are displayed in \Cref{sec:proofs} followed by some technical results in the concluding  part  \Cref{sec:AAL}.

\section{Preliminaries}
\label{sec:prim}
 
\subsection{Lattices}
Unless otherwise stated  $\EIE$ shall denote  a discrete subgroup of the additive group $\TT$ (called also  a lattice on $\TT$) which has finite or countably infinite  number of elements $\abs{\EI}$. In several instances we shall assume  that $\EIE$ has full rank, i.e., for  some non-singular  $l\times l$ real matrix $A$ 
(called base matrix) $\EIE=\{ A x, x\in \TT \},$ where  $x$ denotes a $l\times 1$ vector.  Two base matrices $A,B$ generate the same lattice on $\TT$ if and only if (\ifft)   
$$A= BU,$$ 
where $U$ is an $l\times l$ real matrix with determinant $\pm 1$. 
Therefore, when $\TTT=\R^l$ the volume of the fundamental parallelepiped 
$ \{ A x, x\in [0,1)^l\}$ of the lattice does not depend on the choice of the base matrix and is  given by 
\bqn{ \label{bMat}
	\Delta(\EIE)= \abs{det(A)}>0.
}

\subsection{Some classes of functionals }
\label{some:class}
 Write $\mathcal{H}_\beta\subset \mathcal{H},  \beta\ge 0$ for the class of maps 
(recall the definition of  $\mathcal{H}$   in the Introduction) which are  $\beta$-homogeneous,  i.e.,  $F(c f) = c^\beta F( f)$ for all $f \in \Dset$ and $c>  0$.

Throughout this paper  $\lambda $  denotes the Lebesgue measure on $\R^l$ or the counting measure if the integration is over a countable set. 
Let $g:\TTT\mapsto [0,\IF)$ being locally bounded and  $\lambda$-measurable. When $\TT=\R^l$ it is assumed that $g$ is almost everywhere positive. Write $\mathfrak{G}$ for the class of such $g$'s and write $\Hh_\beta, \beta \ge 0$ for the class of $\beta$-homogeneous maps   
determined for given $\Gamma  \in \mathcal{H}_{\xi_0}, g_i \in \mathfrak{G}, i\le 3$ as
\begin{equation}\label{eier}
	F(f)= \frac{ \Gamma (f) \II_1  (f,g) \ind{\II_2(f,g) \in A}  }{ \II_3  (f,g)  }, \ \  \ \II_i(f, g) = \int_{\TTT} \kappa(B^{-t}f)^{\xi_i}  g_i(t) \lambda(dt), \quad f \in \Dset.
\end{equation}   
Here 
$\xi_i$'s are constants and  $A \subset [0,\IF]$ is such that $F$ is $\beta$-homogeneous.  If $F(f)$ is undefined we set   $F(f)=0$. \\

\subsection{Anchoring \&  involution maps}
Anchoring  maps introduced in \cite{BP} play   a crucial role in the investigation of tail rf's. See also \cite{Hrovje,BP,kulik:soulier:2020, Planic,MR4280158,cissokho2021estimation,kulik2023,chen2023asymptotic} for various results concerning those maps. 
Given  a lattice $\EI$ on $\TT$ we introduce next (positive) shift-involutions and anchoring maps. Below $(\R^k)^*=
 (\R^k \cup \{\IF\}), k\inn$ is the one-point compactification of $\R^k$ and $0$ denotes the origin of $\R^k$ or the zero  of  $\Dset$.
\BD  Let    $\mathcal{J}: \Dset   \to  (\R^d)^* $  \cER{be $\AA/\borel{(\R^d)^*}$-measurable.}  
\begin{enumerate}[{J}1)] 
	\item 
	\label{A1}	For all $j \in \EI,f\in \Dset $ we have 
	$\mathcal{J}(B^j f) = \mathcal{J}( f)+j$;
	\item
\label{A2}  For all $ f\in \Dset$ if $\mathcal{J}(f)=j\in \EI$,
then  $\kappa(B^{-j}f) \ge \min(\kappa(f),1 )$;
	\item
	\label{A3}  For all $f\in \Dset$ if $\mathcal{J}(f)=j\in \EI$,
	 then  $\kappa(B^{-j}f) >0$.
\end{enumerate}	
Suppose that  $\mathcal{J}$  satisfies   \ref{A1}. When  \ref{A2} holds it is referred to as  anchoring. If $\mathcal{J}$ is  0-homogeneous it is called 
a shift-involution and if further \ref{A3} is satisfied it is called a positive shift-involution. 
\ED

 Hereafter  $\prec $ stands for some  total order on $\TT$ which is shift-invariant, i.e., $i\prec j$ implies $i+k \prec j+k$ for all $i,j,k\in \TT$; a canonical  instance  is the lexicographical order. 
 Both  $\inf $ and $\sup $ are defined with respect to  $\prec$. 
 As pointed out in \cite{BP} an   interesting anchoring map is  the first exceedance functional  $\Ife$ defined by 
 $$ \Ife(f)= \inf\Bigl(j\in \EI :  \kappa(B^{-j} f)>1  \Bigr), \quad f \in \Dset,
 $$
 where  $\Ife(f)=  \IF$ if there are infinitely many exceedance on $\{ j \in \EI, j \prec  k_0 \}$ for some 
 $k_0\in \EI$  with  all components positive.  Define further the infargsup map
 $$\Ifasup(f)= \inf \Bigl( j \in \EI:  \kappa(B^{-j}f )= \sup_{i \in \EI} \kappa(B^{-i} f) \Bigr), \ \  f\in \Dset ,$$
  which is  a positive shift-involution \cER{but not} anchoring. The infimum of an empty set is $\IF$. Also if it is not attained at some element of $\EIE$, then the maps defined above are assigned to $ \IF$. \\

\section{Extensions of \eqref{boll}}
\label{sec:shift} 
As in the Introduction consider an  $\alpha$-homogeneous class $\kK$, which is defined 
with respect to a fixed $\alpha>0$ and a given non-negative $\kappa \in\mathcal{H}_\alpha$. 
 Hereafter  $\norm{\vk x}_\IF= \max_{1 \le i \le k} \abs{x_i}, \vk x\inr^k$ and 
 set 
\bqn{ \label{xstar}
	s^{+}=\max(0,s), \quad s^{-}=\max(0,-s), \quad s\inr,  \quad \vk x^{\pm }=( x_1^{\pm}\ldot  x_k^{\pm}),\quad \vk x_{\star}=(\vk x^+, \vk x^-), 
}	
$$\kappa_\IF(f)=\norm{(f(0))_{\star}}_\IF^\alpha, \quad  f\in \Dset.
$$ 
We show next that a simpler condition than \eqref{boll} can be formulated under weak assumptions. 
\BEL \label{parrapT} 
Fix $\vk Z \in \SCD, \kappa \in \mathcal{H}_\alpha$ and suppose that \eqref{ZY} holds  for both $\kappa$ and $\kappa_\IF$. If further	
\bqn{\label{jet} 
	\max( \E*{ \kappa_\IF(B^{t} \vk Z) } , \E*{ \kappa(B^{t} \vk Z) })< \IF, \quad \forall t\in \TTT,
}
 then 
   \eqref{boll} is equivalent with 
\bqn{\label{parrap:E}    \E*{ \max_{1 \le i \le n} \kappa_\IF(B^{t_i} \vk Z /\vk x) }
	&=& 
	\E*{ \max_{1 \le i \le n} \kappa_\IF(B^{t_i} \ZT /\vk x ) }
	< \IF
}
  for all $(t_1 \ldot t_n)\in \diad^n, n\inn,\vk x\in (0,\IF)^d, \ZT\in \kK.$\\
\label{OM1}
\EEL 

\BRM 
	As shown in the proof of \Cref{parrapT},  under the conditions therein $\vk Z$ and $\ZT \in \kK$ are closely related since $\vk Z_{\star}$ and $\ZT_{\star}$ are representers of the same  max-stable rf $\vk X_\star$ and therefore (when $d=1,\alpha=1$) \cite[Eq.\ (3.1)]{MolchanovBE} holds.
\label{remD} 
\ERM

If $F: \Dset \mapsto \R $ is such that $F(\ZT)$ is a rv for all $
 \ZT \in \kK$ and moreover, there exists a sequence of Borel measurable and $0$-homogeneous functions $F_n:\R^n \to [0,\IF]$ such that we have 
 the convergence in probability
 $$
 F_n(\ZT(t_1) \ldot  \ZT(t_n)) \toprob F(\ZT), \quad n\to \IF,
 $$
    with $t_i$'s in $\TT$, then 
 by \eqref{conditionC1} it follows that 
   \eqref{boll} holds also for such general $F$. \\ 
We show next that this is the case for $F$ defined in \eqref{eier} justifying  further the label $\alpha$-homogeneous.

\BT  
For a given $\kK$ \eqref{boll} is equivalent with 
  \bqn{\label{boll3}
	\E{ \kappa(B^h\vk Z) F( \vk Z  )}&=& \E{ \kappa(B^h\ZT) F(\ZT) }, 	\quad 
	\forall F\in \Hh_{0}, \forall  h\in \TT,  {\ZT} \in \kK
} 
and 
\bqn{
	 \label{boll1} 
	\E{ F( \vk Z  )}&=& \E{F( \ZT)}, 	\quad 
	\forall F\in \Hh_{\alpha}, \forall   \ZT \in \kK.
}
\label{lemAHDOOB2}
\ET

\def\JJT{I_\TT}

 \BRM 
\begin{enumerate}[(i)]   \label{borsh23}
\item We interpret 	$\IF \cdot 0$ and  $0/0$ as 0. This is the case for instance if the  maps $F$ above are products or  ratios with one component being the indicator function; 
\item \label{remD1}In view of the above result also $\kappa \in \Hh_\alpha$ can be considered. Particular choices of interest (as suggested from one referee) are 
$$ 
\kappa(f)= \max_{t\in T_0} \norm{f(t)}^\alpha,
$$
with  $T_0 \subset \TT$ having countable number of elements.

\item So far the assumption \eqref{ZY} has been used in the derivation of \eqref{boll1} only. If  $F \in \Hh_{\alpha}$ is such that 
$$
\E*{F(\ZT) \ind*{\sup_{t\in \TT} \kappa(B^{-t} \ZT)=0}}=0, \quad \forall \ZT \in \kK,
$$
  then \eqref{boll1} still holds without \eqref{ZY}. When $\TT$ has countable number of elements instead of \eqref{ZY} we can assume $F(0)=0$. 
\item By  \eqref{boll3} with pdf $p_N(t)>0, t\in \TT$ we have (set $	\JJT(f)= \int_\TT \kappa(B^{-t}f) p_N(t) \lambda(dt)$)
$$\E*{F(\vk Z) \JJT(\vk Z)}=\E*{ F(\ZT)\JJT(\vk Z)}, \quad \forall F\in \Hh_0, \ZT\in \kk.
$$
 Hence from \eqref{ZY} and applying \Cref{lemGK} taking therein  $d=1,\mathcal{A}=\R^l,\vk U=\kappa(B^{-t} \ZT), \gamma=z=0$ and $g_1(s)=s^\alpha,s\ge 0, g_2(t)=p_N(t),t\in \TT$ we obtain $\JJT(\vk Z)>0$ a.s.\ and therefore     
\bqn{\label{prosimo}
	\E{F(\vk Z)} =0 \quad \implies  \E{F(\ZT)}=0, \quad \forall \ZT\in \kK.
}
If $\kK$ consist of rf's with \cadlag\ sample paths, for any measurable cone $C \subset \Dset$ (see \cite{Domb14} for the definition), 
 by \eqref{prosimo} 
\bqn{\label{prosimo2}
	\pk{ \vk Z \in C} =1 \quad \implies  \pk{\ZT \in C}=1, \quad \forall \ZT\in \kK
}
and thus we retrieve \cite[Prop 1]{Domb14}.
\end{enumerate}
\ERM

 Borrowing the idea of \cite[Thm 4.2]{WangStoev}, the next application of \Cref{lemAHDOOB2} shows how to split $\kK$ into shift-invariant classes.
\begin{example}[Splitting of $\kK$]
	Given a $\kK$ suppose that $\vk Z= \vk Z_1+ \vk Z_2,$ where $\vk Z_i$'s belong to  $\SCD$ are such that  both $\vk Z_1, \vk Z_2$ satisfy    \eqref{ZY}, \eqref{lBo} and \eqref{conditionC1}. Hence we can define 
	$ \mathbb{K}_\alpha[\vk Z_i], i=1,2 $.  
	It follows easily that 
	under the following condition (which is borrowed from \cite[Eq. (13.2.12)]{kulik:soulier:2020})
	\bqn{ \label{Gtw}
		\pk*{  \sup_{t\in \TT} \kappa(B^{-t} \vk Z_1) =0, \sup_{t\in \TT} \kappa(B^{-t} \vk Z_2) =0}=0
		}
we obtain for all $\ZT_i \in  \mathbb{K}_\alpha[\vk Z_{i} ], i=1,2$ defined in the same probability space and satisfying \eqref{Gtw}  that
	$$  \ZT_1+  \ZT_2  \in \kK.
	$$
If $\vk Z$ has non-negative components and $\kappa(f)= \norm{f(0)}^\alpha$
 (maximum is taken component-wise below) 
\bqn{
	\label{shuto} \vk X(t)= \max( \vk X_1(t), \vk X_2(t)), \quad t\in \TT
}
is a max-stable rf as in \eqref{eq1} with representer $\vk Z$, if 
$\vk X_k$ has the same de Haan representation \eqref{eq1} with representer $\vk Z_k, k=1,2$. 
In view of \eqref{kode} and \eqref{Gtw}  it follows that $\vk X_1$ and $\vk X_2$ are independent. This splitting property is well-known, see \cite[Lem 5]{dom2016}, 
\cite[Lem 13.2.8]{kulik:soulier:2020} and \cite[Thm 4.2]{WangStoev}. 
\label{examp:2bc}
\end{example}
\begin{example}[Cone splitting of $\kK$] 
Let   $C_i \subset \Dset, i \in I=\{1 \ldot L\}$ with $L$ a positive integer or equal to infinity such that $H_i(c\ZT)=\ind{ c\ZT\in C_i}$ is a rv for all $c>0, \ZT \in \kK$ and  
$$\pk{H_i(c\ZT ) =H_i(\ZT )}=1,  \quad \forall c>0, \forall i\in I .
$$ 
Suppose that $p_i=\pk{ \vk Z \in C_i} \in (0,1), i \in I$ and define  $\vk Z_i$ to be equal in law with $\vk Z$ conditioned on $\vk Z \in C_i$. Write $N$ for a rv taking values in $I $ such that $\pk{N=i}=p_i, i=1 \ldot L$ and let $$\ZT_i \in \mathscr{C}_\kappa[\vk Z_i], \quad i\in I $$
 be  defined in the same probability space being independent of $N$. 
If further $ \pk{\sum_{n\in I} \ind{ \vk Z\in C_n}=1}=1$, then it follows from the choice of $C_i$'s and \eqref{boll1} 
that 
\bqn{
	\label{eq:2bd} \ZT_N  \in \kK.
}
\label{examp:2bd}
\end{example}
\section{Characterisation of shift-invariance}
\label{sec:char:shift} 
In order to deal with the shift-invariance property, we shall assume that $\TT$ is an additive group, whenever this property is mentioned. 

By definition, $\kK$ is shift-invariant if $B^h \vk Z$ belongs to $\kK$ for all $h\in \TT$. 
Consequently,  by \eqref{tcfN0} and \eqref{boll3} $\kK$ is shift-invariant \ifft
\bqn{  \E*{ \kappa( B^{-h}\vk Z)  F( \vk Z )} &=&  	\E*{ \kappa( \ZT)F(B^h \ZT) } , \label{tcfN00} 
	\ \  \forall F\in \Hh_{0},\forall h\in \TT,  \forall \ZT\in \kK,
} 
 which  in view of \eqref{boll1} is also equivalent with 
\bqn{ \label{tcfN} 
	\E{ F({\vk Z})} = \E{F( B^h  \widetilde{\vk Z})}, \quad \forall F\in \Hh_{\alpha} , \forall h\in \TT, \forall \widetilde{\vk Z}\in \kK.
} 
In the sequel  a shift-invariant $\kK$   will be simply denoted by $\kk$. Assume next without loss of generality that 
\begin{equation}\label{pobien}
	\E{ \kappa({\vk Z }) }= 1.
\end{equation}  
Hence, for a given   $\kk$ 
\bqn{\label{eqYz}
1 &=&\E*{\kappa(B^h \ZT) }, \quad \forall h \in \TT,  
	\forall \ZT \in \kk.
}	
\BRM
Under the assumptions of \Cref{OM1}, if  further $\kK$ is shift-invariant, then in view of the proof of \Cref{OM1} $B^h \vk Z_{\star}, h\in \TT$ 
and $\vk Z_{\star}$ are two representers of  the shift-invariant max-stable rf $\vk X_{\star}$. These representers are directly related, see 
\cite[p.1219, case $d=1$]{MolchanovBE} and \cite[Thm 1, case $d\ge 1$]{MR2453345}. Consequently, also $B^h Z$ and $Z$ are related, recall the definition of $x_\star$ in \eqref{xstar}.
\ERM  
Clearly, by \eqref{tcfN} all $\kk$'s  are closed under random shifting, i.e.,   
\bqn{\label{provedKab}\vk Z_N=\ B^N \vk Z \in \kk
}
for all  $\TT$-valued rv's $N$ defined in $(\Omega, \mathscr{F}, {\mathbb{P}})$ being independent of $\vk Z$.
 For  special $\vk Z$ included in \Cref{examp:5} the claim in \eqref{provedKab} 
 is proved in \cite[Thm 2]{KabExt}.
    
\begin{example}[\Shiftinvariance of Brown-Resnick  $\kK$]
	Consider the settings of \Cref{examp:5} \cER{and suppose for   simplicity that $\alpha =1$}.  
		 In view of \cite[Thm 22]{lorenzo} the pseudo-variogram matrix-valued function $\gamma$ is conditionally negative definite. Moreover, examples of such vector-valued Gaussian rf's $\vk V$  do exist. Suppose next    that $\gamma(s,t),s,t\in \TT$ 
		 depends only on $t-s$ for all $s,t \in \TT$. For given positive integers $i,k \le d$  by 
\Cref{lemGauss}	the tilted law of $ W_{j}(t)-W_k(h)$ with respect to  $e^{ W_k(h) }$
	is the same as that of (recall  $W_{j}(t)= V_{i}(t)- Var(V_{j}(t))/2$)
\bqn{ W_{j}(t)-W_k(h)+ Cov(V_j(t),V_k(h))-Cov(V_k(h),V_k(h))= 
	 V_{j}(t)-V_k(h)- \gamma_{jk}(t,h).
	 \label{cohen}
}	 
	By the assumption on $\gamma$ the latter has the same law as  $ V_{j}(t-h)-V_k(0)- \gamma_{jk}(t-h,0)$ 
	and moreover \eqref{cohen} holds jointly for all $1 \le j\le d $. Consequently, since $\xi_i$'s are independent of $\vk V$, 
 	for all $F \in \Hh_0$ using the 0-homogeneity of $F$ in the derivation of the  first and the third  equality below 
	\bqny{  
		d\E{ \norm{\vk Z(h)}^\alpha_\alpha  F(\vk Z)}& =&
			\sum_{k=1}^d 	 \E{ e^{   W_k(h)} F(\vk Z\cER{/Z_k(h)})) } \\
			&= &
			\sum_{k=1}^d 	 \E{ e^{    W_k(0)} F(B^h\vk Z\cER{/Z_k(0)}) }\\
			&=& 
				 \E*{ \sum_{k=1}^d \abs{Z_k(0)}  F(B^h\vk Z) }\\
			&=& 
			d\E{\norm{\vk Z(0)}^\alpha _\alpha F(B^h \vk Z) }, \quad \forall h\in \TT.
		} 
		Hence $\kK$ is shift-invariant, which  for $d=1$ is a direct consequence of \cite{kab2009}. 
	In view of \eqref{cohen} and the derivation above the law of the local rf $\vk \Theta$ depends only on $\gamma$  (even without imposing further restrictions on $\gamma(s,t)$), which is shown for $d=1$ already in \cite{kab2009}.\\ 
Note  that when $\xi_i$'s are deterministic, the \shiftinvariance of $\kK$ is equivalent with the stationarity of the max-stable   rf $\vk X$ defined in the Introduction. Under the settings of this example, \cite[Lem 4.2]{KumeE} shows the stationarity of $\vk X$. 
		\label{examp:5B}
	\end{example}	 
\begin{example}[$\kk$'s defined from pdf's] Consider $\xi_i$'s as in \Cref{examp:5} and let $\varphi_i,1 \le i \le d$ be pdf's of $k$-dimensional centered Gaussian random vectors and set $\vk \varphi=( \xi_1\varphi_1 \ldot\xi_d\varphi_d ) $. Assuming that $N$ is independent of $\xi_i$'s,   
	setting 
	$$
	\vk Z_N(t)= [p_N(t)]^{-1/\alpha} \vk \varphi (t-N), \quad t\in \TT
	$$
	 and taking $\kappa $ as in \Cref{examp:5}, it follows that $\vk Z_N$  satisfies 
	\eqref{ZY}, \eqref{boll} and \eqref{conditionC1} and the corresponding $\kK$ is shift-invariant.
\end{example}

\BD Write $\pclas $   for the class of functions 
$$p_N: \diad \mapsto (0,\IF), \quad \intDD p_N(t)\lambda(dt) < \IF,$$
where $\lambda$ is the counting measure on countable dense set $\diad$ (defined in the Introduction). We define similarly $\pclasT$ when $\TT=\R^l$ substituting $\diad$ by $\TT$ and 
taking $\lambda$ to be the Lebesgues measure on $\TT$  assuming further that  $p_N$ is almost everywhere positive and locally bounded.   
\ED

We state next   equivalent characterisations of shift-invariant $\kK$'s.  
Write $\vk Z_\diad(t)=\vk Z(t), t\in \diad$ for the restriction of $\vk Z$ on $\diad$ and denote the corresponding $\alpha$-homogeneous class of rf's  by 
$\mathscr{C}_\kappa[\vk Z_\diad]$. 
In the following $N$ is  independent of any other random element, being  a $\diad$-valued or $\TT$-valued rv with density $p_N \in \pclas$ or $p_N \in \pclasT$, respectively. \\
In order to define the integral
$$ I_\TT(\vk Z)= \intT \kappa(B^{-r}\vk Z) p_N(r) \lambda(dr)$$
when $\TT=\R^l$ we need the  joint measurability of $\kappa(B^{-t}\vk Z), t\in \TT$, which is a direct consequence of the stochastic continuity assumption. The latter is also important for the positivity of $I_\TT(\vk Z)$ when applying \Cref{lemGK} below.  

\BEL\label{shpejt}    
 $\kK$ is shift-invariant iff one of the following holds:	 
\begin{enumerate}[(i)] 
	\item \label{item:shpejt:0} $ \mathscr{C}_\kappa[\vk Z_\diad]$ is shift-invariant;  
	\item \label{item:shpejt:00} \eqref{tcfN00} 
	is satisfied for all $h\in \diad$;
\item \label{item:shpejt:2} for all rv's $N$ with $p_N \in \pclas$  the following rf's
$$ \vk Z_N(t) =  \Bigl(\frac{\tX}{  I_\diad(B^N \vk Z) }\Bigr)^{1/\alpha}
 B^N \vk Z(t) 
 $$
 satisfy \eqref{boll3};
\item \label{item:shpejt:4}  \Cref{item:shpejt:2} holds with $\diad$ substituted by $\TT$ for all $p_N \in \pclasT$.
\end{enumerate}
\EEL

\BD   $F\in \wHO$ (and similarly for other maps) is shift-invariant with respect to  $\EI\subset \TTT$ if 
$$F(B^hf)= F(f), \quad\ \forall f\in \Dset,  \forall h\in \EI.$$
\ED 

Below shift-invariant maps are considered with respect to $\EI=\diad$.  Therefore   we shall simply write  shift-invariant maps without specifying $\EI$.
\begin{example}[Shift-invariant maps \& new shift-invariant classes of rf's]
		Consider \  $\kk$ and 
		a \shiftinvariant\  non-negative map $\Gamma\in  \Hh_0$. 
		If 
		$\E{ \kappa(\vk Z) \Gamma(\vk Z)}\in (0,\IF)$ by  \eqref{tcfN}  
		$\mathbb{K}_\alpha[ \hat {\vk Z}]$ is also \shiftinvariant\ with $\hat{\vk Z }$ equal $\vk Z  \Gamma(\vk Z)$ conditioned on $\Gamma(\vk Z)>0 $. 
		\label{examp:2b}
	\end{example}
\begin{example}[Cone splitting of $\kk$] 
		Reconsider \Cref{examp:2bd}  assuming additionally that $H_i$'s therein are shift-invariant. In view of \eqref{tcfN00} and \Cref{shpejt}, \Cref{item:shpejt:0} we have that $\mathscr{C}_\kappa[\vk Z_i]$'s  are shift-invariant. Moreover, taking 
		independent $\ZT_i$'s from $\mathscr{C}_\kappa[\vk Z_i]$'s  defined in the same probability space as $N$ 
		being also independent of $N$, then we have 
		\bqn{\label{eq:2bde} \ZT_N  \in \kk.
		}
		\label{examp:2bde}
	\end{example}

\section{Tail and spectral tail RF's }
\label{secTailSpectral}
Recall that  given a $\kk$ we define a local rf    $\widetilde{\vkT}$ as the  rf $\ZT/\kappa(\ZT )^{1/\alpha} $ under $\widehat{\mathbb{P}}$ given by 
	\eqref{minist}, where $\ZT \in \kk$.  		
Note that our assumption $(\Omega, \mathcal{F}, {\mathbb{P}})$  is complete yields that also 
		$(\Omega, \mathcal{F},\widehat{\mathbb{P}})$ is complete.  For notational simplicity we    write  $ {\mathbb{P}}$ instead of $\widehat{\mathbb{P}}$. In view of \eqref{boll1} and \eqref{eqYz}     for local rf's $\vkT$ and $\widetilde{\vkT}$
\bqn{\label{eqDo20}
\E*{\kappa( B^{-h} \widetilde{\vkT} )  \Gamma(\widetilde{\vkT} )}=	\E*{\kappa( B^{-h} \vkT) \Gamma( \vkT)}
&=& \E*{ \mathbb{I}( \kappa(B^h   \vkT) \not =0) \Gamma( B^{h} {\vkT}) }, 
\quad  \forall \Gamma\in  \wHO,  \forall h\in \TT.
} 
Consequently, the fidi's of $\vkT$ and $\widetilde{\vkT}$ agree. Hence in view of \Cref{lemAHDOOB} if $F\in \mathcal{H}_\beta, \beta \ge 0$, the law of $F(\vkT)$ is the same as that of $F(\widetilde{\vkT})$.\\
For the special case $\kappa(f)= \norm{f(0)}^\alpha,  f\in \Dset$ the identity \eqref{eqDo20} is equivalent with the time-change formula discovered in \cite{BojanS}. Therein 
$\vk \Theta$ is the so-called spectral tail rf of a regularly varying time series. The time-change formula for a  radial function $\kappa$ appears in \cite[Thm 6.1]{Segers}. 
\BRM \label{rem:unique}
If $\pk{\kappa(\vk Z)>0}=1$, then since we assume \eqref{pobien}, from \eqref{eqDo20} we conclude that  $\vkT \in \kK$. For this case a similar formula to \eqref{eqDo20} referred to as the non-singular group action property is shown in \cite[Thm S.10]{klemus}.
\ERM 
\def\vkT{\vk \Theta}
  Motivated by the spectral tail rf's appearing in the context of time series, we define a  corresponding rf without any reference to a particular $\kk$.

\BD $\vk \Theta \in \SCD$  
is called a spectral tail rf, if
\begin{enumerate}[(i)]
	\item \label{defA:2} $\pk{\kappa(\vkT)=1}=1$;
	\item \label{defA:1} for all   $\Gamma\in  \wHO,  h\in \TT$ the identity \eqref{eqDo20} holds; 
	\item \label{defA:3}  
	for all $p_N\in \pclas$ (recall $I_\diad(\vkT)=\intDD \kappa(B^{-r} \vkT) p_N(r) \lambda(dr)$)
\bqn{ 
	\E*{ \frac{ \sup_{s \in [-c,c]^l\cap \diad} \kappa(B^{N-s} \vkT )}{ 
			I_\diad(B^N \vkT)}}< \IF, \quad \forall c>0. 
	\label{justen}
}  
\end{enumerate}    
\label{defA}
\ED 
\BRM 
\begin{enumerate}[(i)] 
\item \cER{If $\TT=\Z^l$, then $\diad=\TT$ and hence \eqref{justen} is satisfied since $[-c,c]^l\cap \diad$ has only finite number of elements and for all $s\in \TT$ by \eqref{eqDo20} and $\pk{\kappa(\vkT)=1}=1$
\bqny{ 
	\E*{ \frac{ \kappa(  B^{N-s}\vkT)  }{		{ I_\diad( B^N\vkT) }}} 
	&=& 
 \intDD\E*{ \frac{ \kappa( B^{s-t}\vkT) }{  I_\diad( B^s\vkT) }}  p_N(t) \lambda(dt)  
= \E*{ \frac{ \intDD \kappa( B^{s-t}\vkT) p_N(t) \lambda(dt)}{ 		{ I_\diad( B^s\vkT) }
}}  =1.
	\label{justen2}
}  
Note that $I_\diad(\vkT) \in [p_N(0),\IF)$ a.s.\ follows by the definition of $p_N$ and  $$\E{\kappa( B^h \vkT)}=\pk{ \kappa(B^h \vkT) \not=0}\le 1, \quad h\in \TT.$$
\cER{Moreover, $I_\diad(\vkT) $ and in general  $F(\vkT), \forall F\in \Hh_\beta, \beta\ge 0$ is well-defined as shown in  the proof of \Cref{lemAHDOOB2}}.
\item An equivalent condition to \eqref{justen} is given in \eqref{ey} below, see for a more restrictive setup  \cite{PH2020,MartinE}.
}
\end{enumerate}
\ERM
 
\def\ZY{\vk Z^*}

\BEL   Let $\vkT$ be the local rf of $\kK$.
If   $\vkT$ is  a spectral tail rf   and for all rv's $N$  independent of $\vkT$ such that 
\bqn{  \label{insintonia}
\vk Z_N(t) =  \frac{1}{  \bigl(I_\diad( B^{N}  \vkT) \bigr)^{1/\alpha}} B^N \vkT(t) \in \kK,
}
$p_N \in \pclas$  $p_N \in \pclasT $  (with $\diad$ substituted by $\TT$), then 
	$\kK$ is shift-invariant. Conversely, if  $\kK$ is  shift-invariant, 
	then its local rf $\vk \Theta$ is   a spectral tail rf  and moreover \eqref{justen} holds also for $\diad$ substituted by $\TT$  and all $p_N \in \pclasT$.
		\label{theoA}
\EEL 
	
\BRM
\label{remKA} 
 	Given a spectral tail rf $\vkT$, the  class of rf's $\mathscr{C}_\kappa(\vkT)$ of elements from $\SCD$ that includes all $\vk Z_N$ determined by 
\eqref{insintonia} for both $p_N \in \pclas$ or $p_N \in \pclasT$ and all $\ZT\in \SCD  $ such that \eqref{boll} holds with $\vk Z=\vk Z_N$ is  $\alpha$-homogeneous and shift-invariant. 
\ERM
Hereafter the $\alpha$-Pareto rv  $ R$ with survival function $t^{-\alpha}, t\ge 1$  is assumed to be independent of any other random element. For $\vk Y= R \vk \Theta$   
using \cite[Lem 3.16]{MartinE} or as in \cite{Hrovje,BP}, we obtain
\bqn{  
	\E{ \Gamma( xB^h \vk Y) \ind{ \kappa(x B^h \vk Y) >1 }} &=& 
	x^\alpha \E{  \Gamma(  \vk Y) \ind{ \kappa(B^{-h}\vk Y/x) >1 }   }, \quad  \forall \Gamma \in \mathcal{H}, \forall h\in \TT, \forall x>0.
	\label{tYY}
}  
\BD $\vk Y \in \SCD   $ is called a   tail rf,  if
\begin{enumerate}[(i)]
	\item \label{defB:2} $\pk{\kappa(\vk Y)> 1}=1$;
	\item \label{defB:1} \eqref{tYY} holds for all $x>0$;
	\item \label{defB:3} Eq.\ \eqref{justen} holds for all $p_N \in \pclas$ with $\vkT$ substituted by $\vk Y $.   
\end{enumerate}      
\label{defB}
\ED 
Tail and spectral tail rf's are crucial for the study of  stationary regularly varying time series, see \cite{BojanS, Hrovje,SegersEx17, klem,WS,kulik:soulier:2020}. 
\BD $\wwH $ consists of all maps $F$ defined by 
$$F=F_1F_2, \ \ F_1\in\wH_\beta ,\ F_2 \in \mathcal{H}, \quad \forall \beta\ge 0 .$$
\label{defH}
\ED   

\BRM \label{remalso} 
\begin{enumerate}[(i)]	\item
	When $\kappa(f) =\norm{f(0)}^\alpha$, then \eqref{tYY}   is shown in \cite{Hrovje, BP} assuming  $\TT= \Z^l$ and in \cite{PH2020,MartinE} considering the \cadlag\ case $\TT=\R$ and $\TTT=\R^l$, respectively.
	\item \label{tYYY}
	 If  $\vk \Theta$ is the local  rf   of some shift-invariant $\kk$, 
	then $\vk Y= R\vk \Theta$ is  a tail rf. Moreover, as in the proof of \cite[Appendix A]{Hrovje} 
	  \eqref{eqDo20} implies that \eqref{tYY} holds for all $\Gamma \in \wwH$.
	
	\item \label{remalso:2} $\Gamma$ in \eqref{tYY} can be of the form 
	$\IF \cdot 0$ or  $0/0$ and both interpreted  as 0.
\end{enumerate}  
\ERM   

The next result  is motivated by  \cite[Thm 3.2]{WS} and \cite{PH2020} which shows \eqref{ey} for  $l=1, \tau=0$ in the  \cadlag\ case. The values of $\tau$ considered below are stated in the next condition: 
\begin{condition} \label{CR2} 
	 $\tau \inr$ is such that  
		\bqn{ \sup_{t\in [-c,c]^l \cap \TT} \E{\kappa( B^{-t}\vk \Theta)^\tau} < \IF, \quad \forall c>0. 
			}
\end{condition}

\BT If $\vk Y$ is a tail rf, then $R=\kappa( \vk Y)^{1/\alpha}$  is an $\alpha$-Pareto rv  being independent of 
$\vk \Theta= \vk Y/ R$. Further   $\vk \Theta$ is a spectral tail rf  satisfying 
\bqn{  \ \ \ \int_0^\IF \E{ \Gamma(z \vk \Theta) 
		\ind{  \kappa(z B^{-h}\vkT)> 1}} \drz=
	\E{  \Gamma(  B^h\vk Y) \ind{  \kappa( B^h \vk Y) \not=0 }   }, \ \forall \Gamma \in \wwH , \forall h\in \TT.
	\label{tYYGen}
}  
Moreover, for $ \tau$ satisfying  \Cref{CR2} there exists a $\kk$ with local rf $\vk \Theta $ such that
\begin{equation}\label{ey}
	\int_{t\in [-c,c]^l\cap \TTT} \E*{\frac{1 }{\int_{s \in [-c,c]^l \cap \TTT} \kappa( B^{t-s}\vk \Theta)^\tau \ind{ \kappa(B^{t-s}\vk Y) >1}\lambda(ds)}}\lambda(dt)< \IF, \quad \forall c>0.  
\end{equation} 
Conversely, if  $\vk \Theta  \in \SCD $ satisfies Definition \ref{defA}, \Cref{defA:2}  and 
both \eqref{tYYGen}, \cER{\eqref{ey} hold with $\vk Y=R\vk \Theta$, then $\vk Y$  is a tail rf.}
\label{altFiorentinischeTrad}
\ET

\def\kkat{\mathbb K^*_\alpha(\vk \Theta) } 
\def\KE{\mathcal K_\alpha } 
\def\setA{\mathcal{A}}

\section{New $\alpha$-homogeneous classes \& Universal Maps}
\label{sec:uni}
In this section  $\EI \subset \TT$ is a discrete     subgroup of the additive group $\TT$  (or simply a lattice). Given a $\kk$, in the light of our results above, the following  approaches lead to new classes of  shift-invariant $\alpha$-homogeneous  rf's:  
\begin{enumerate}[(i)]
	\item $\mathscr{C}_{\kappa}[\vk Z_{ \EI} ]$ is shift-invariant, where $\vk Z_{ \EI}(t)= \vk Z(t),t\in \EI$;
	\item Splitting with respect to  $C_i$'s as in \Cref{examp:2bde} leads to new shift-invariant $\mathscr{C}_{\kappa}[\vk Z_i]$'s such that further \eqref{eq:2bde} holds;
	\item Given shift-invariant  $\mathcal{C}_{\kappa}[\vk Z_i], i=1 \ldot L$ with $L$ a positive integer or equal to infinity, then taking $\ZT_i \in \mathcal{C}_{\kappa} [\vk Z_i],i=1 \ldot L$ in the same probability space  independent of a rv $N$ defined in the same probability space, then $\mathscr{C}_{\kappa}[\ZT_N]$ is also shift-invariant;
	\item Utilising $N$ and $\vk Z_N$ as in \Cref{shpejt} and then defining the corresponding $\mathscr{C}_{\kappa}[\ZT_N]$ by keeping only $\ZT$'s that satisfy \eqref{ZY}.
\end{enumerate}

The last approach to construct shift-invariant $\kk$'s relies on the fact that by choosing  a given $p_N$, 
it is possible to show that the corresponding integral functional of $\vk Z$  is a.s.\ positive and finite. 

It is possible to utilise functionals that do not depend on the choice of some given $p_N$. A particular class of such functionals, referred below as universal maps, play a crucial role in the construction of random shift-representations which are studies in the subsequent contribution \cite{EHPara}.  \\ Hereafter all event inclusions/ equalities are modulo null sets with respect to the corresponding probability measure. Let  $\cEp{V\subset \TT}$  be such that $\lambda(V)>0$. For  given $f: \TTT\mapsto \cED{\R^l}$    set 
	$$ \mathcal{S}_{V}(f) = \int_{V} \kappa(B^{-t} f)  \lambda(dt), \quad   
	 \mathcal{B}_{V, \tau}(f)=
	\int_{V}  \kappa( B^{-t} f) ^{\tau}\ind{ \kappa(B^{-t} f) > 1 }\lambda(dt).
$$
When the above are not properly defined we set them equal to 0. If $V=\TTT$ we write  $\mathcal{S}$ instead of $\mathcal{S}_V$.

\def\SLI{\mathcal{S}(\vk Q_\IF)}

In the rest of this section 
 $\EI$ has countably infinite number of elements. A shift-involution   will be denoted   by $\mj_1$, whereas  $\mj_2$ and $\mj_3$  denote 
a positive shift-involution and an anchoring map, respectively (with respect  to   $\EI$). 
Given a tail rf $\vk Y$ a weak restriction introduced in \cite{HBernulli} that  we shall  impose next   is   
\bqn{  \pk{ \mathcal{S}_\EIE(\vk Z) < \IF}= \pk{\mathcal{S}_\EIE(\vk Z) < \IF, \mj_1(\vk Z )  \in \EI}
	\label{peres1}
}	
and 
\bqn{  \pk{ \mathcal{S}_\EIE(\vk Y) < \IF}= \pk{\mathcal{S}_\EIE(\vk Y) < \IF, \mj_k(\vk Y )  \in \EI}, \quad k=2
	\label{peres2}
}
and further $\mj_3$ satisfies \eqref{peres2} with $\mathcal{S}_{\EIE}(\vk Y)$ substituted by $\mathcal{B}_{\EIE,0}(\vk Y)$. These two conditions are satisfied in particular  if $\mj_i=\Ifasup, i=1,2,3$ and $\mj_3= \Ife$.    
Below we set $\norm{t}_*=\sum_{i=1}^l \abs{t_i}, t=(t_1 \ldot t_l)\in \R^l$.\\
\BD Denote by $\clasU$ the the class of maps  $F: \Dset \to (\R^k)^*, k\inn$   such that for all $F_1,F_2\in \clasU$ there is a measurable subset  $\setA$ of $(\R^k)^*$ so that   for all   tail rf's $\vk Y$ we have 
$ \{ F_1(\vk Y) \in \setA \} =  \{ F_2(\vk Y) \in \setA \}$. We call $\clasU$ the class of universal maps.  
\ED 
 The next lemma shows in particular that $\clasU$ contains both $\mathcal{S}_{\EI}(\cdot)$ and $\mathcal{B}_{\EI,\tau}(\cdot)$.  
\BEL  
Let $\EI$ be a discrete subgroup of the additive group $\TT$  and suppose that $\abs{\EI}=\IF$. 
If $\mj_1$ satisfies
\eqref{peres1}, $\mj_2$ satisfies \eqref{peres2} and $\mj_3$ satisfies \eqref{peres2} with  $\mathcal{S}_{\EIE}(\vk Y)$ being substituted by $\mathcal{B}_{\EIE,0}(\vk Y)$,  then   
\bqn{ 
	\quad \quad \mbox{    } 
\{  \SZq < \IF \} =	\{ \mj_1(\vk Z) \in  \EI \}= 	\Bigl \{ \lim_{ \norm{t}_*  \to \IF, t\in \EI} \kappa(B^{-t} \vk Z) =0   \Bigr \}, 	\label{ideT00}\\   
\ 	\{ \mathcal{I}(\vk Y) \in  \EI\}= \{ \STYq < \IF \} = \{\JJqx < \IF \}= \Bigl \{ \lim_{ \norm{t}_* \to \IF, t\in \EI} \kappa( B^{-t} \vk Y) =0   \Bigr \} , \label{ideT0}
	\label{ideTE}
}	
with $\mathcal{I}$ equal to $\mj_2$ or $\mj_3$ \cER{and all $\tau\inr$.}
\label{domodosola} 
\EEL

\def\BBY{\mathcal{B}(\vk Y)}  
\def\BBYt{\mathcal{B}_{\TTT,\tau}(\vk Y)}  
\def\MM{\mathcal{M}(\vk Y)}
\def\MME{\mathcal{M}_\EIE(\vk Y)}
  
Since $\pk{\kappa(\Theta)=1}=1$ and $0 \in \EI$, then clearly $\pk{\STq>0}= \pk{\JJqx>0  }=1$. In view of 
\Cref{lemAH} it follows further that $\pk{\SZq>0}=1$. When $\EI$ is substituted by $\TT$, these conclusions are not trivial.\\
We present below the main result of this section. 
 
\def\SZT{\mathcal{S}(\ZT)}

\BT Under the assumptions of \Cref{domodosola}, if  $\tau$ satisfies   \Cref{CR2}, then we have 
\bqn{ 
		\quad 	\pk{\SZ>0}=	\pk{\ST >0}=\pk{\BBYt >0}=1 . 
		\label{am0}
}	
Moreover, if $\EI $ is a full rank lattice on $\TT$, then we have    
\bqn{ \label{am10}
\qquad	
	\{  \ST < \IF \} 
	 & =&  	\Bigl \{ \lim_{ \norm{t}_* \to \IF} \kappa( B^{-t}\vk \Theta ) =0   \Bigr \}
	 =   \Bigl\{ \intT \sup_{t \in  K \cap \TTT  } \kappa( B^{s-t}  \vk \Theta)\lambda(ds) < \IF \Bigr\}	=\{\JJqx< \IF  \}\\ 
	 & =& \Bigl\{ \BBYt < \IF \Bigr\} 
	=	 \Bigl\{ \STYq < \IF \Bigr\}
	=	\{ \mathcal{I}_1(\vk Y ) \in \EI \}=\{ \mathcal{I}_2(\vk Y ) \in \EI \}  		 \label{am12}
}	
and  for all $\ZT\in \kk$ and all compact sets $K\subset \R^l$
\bqn{ \label{am20}
	\{  \SZT < \IF \} 
	& =&  
	\Bigl \{ \lim_{ \norm{t}_* \to \IF} \kappa(B^{-t}\ZT)=0   \Bigr \}  	=
	 \Bigl\{ \intT \sup_{t \in  K \cap \TTT  } \kappa(B^{s-t}\ZT) \lambda(ds) < \IF \Bigr\}\\ 
	& =& 		\{ \mathcal{I}_1(\ZT ) \in \EI \}= \{\mathcal{S}_{\EIE}(\ZT)< \IF  \}
	 \label{am20c}
 }	
 hold with $\mathcal{I}_1=\Ifasup$ and  $\mathcal{I}_2=\Ife$.
\label{kabA}
\ET

\BRM 
\begin{enumerate}[(i)]
	\item \Cref{kabA} extends several claims derived for particular choices of $\TT, \EI, \tau=0$ and $\kappa(f)=\norm{f(0)}^\alpha$  in 
 \cite{dom2016,kab2009a, WangStoev, BP, WS, HBernulli,klem,Hrovje, kulik:soulier:2020, PH2020}. In case of max-stable or symmetric $\alpha$-stable rf's, the results in the literature that we extended above concern also several criteria for pure dissipativity/conservativity. Similar characterisations also appear in the analysis of regularly varying rf's. 
\item \cER{The importance of lattices in the study of dissipativity of symmetric $\alpha$-stable rf's is know from \cite{Roy1}.}
\item \cER{The class $\clasU$ is non-empty, since all functionals in \Cref{domodosola} and \Cref{kabA}   belong to it (under the conditions therein).} 
 \label{lemXX}
\end{enumerate}
\ERM 
     
\def\KK{ \mathscr{K}}

 	\section{Discussions}
 	\label{sec:impli}
 We consider in this section some implications for max-stable rf's and 
calculation of maximal indices. 
 
\subsection{Max-stable rf's}
 \label{sec:maxstable}
 	A given $\kK$ with $\vk Z$ having non-negative components and local rf $\vk \Theta$  is closely related to max-stable rf's as shown in \cite{KumeE}.  For such $Z$, define
 	 a max-stable rf $ \vk X(t),t\in \TT$ as in \eqref{eq1}. \\
 The tractability of $\vk X$ is related to expression \eqref{kode}, which can be rewritten   as 
 \bqny{ 
 	\lefteqn{	- \ln \pk{ \vk X(t_i) \le \vk x_i, 1 \le i\le n}}\\ 
 	&=&	  \sum_{1\le l \le n} \E*{   \norm{ \vk Z(t_l)/\vk x_{l}}^\alpha _\IF
 		\mathbb{I}( infargmax_{  1 \le j \le n }   \norm{ \vk Z(t_j)/\vk x_{j}}^\alpha_\IF =l  )  }
 	} 
 	for all $t_i \in \TT, 
 	\vk x_i\in (0,\IF)^d, i=1 \ldot n$. 
Moreover,  if  $\kappa$ is such that a.s.\ 
 	\bqn{\label{re}
 	\{\kappa(\vk Z) \not=0\}= \{\norm{\vk Z(0)}\not =0\}
 }
 and $\kK$ is shift-invariant, then  another expression for \eqref{kode} is shown in \cite{Htilt}. Namely,  as therein for general $\vk \Theta$ defined with respect to this $\kappa$, we obtain by \eqref{tcfN} 
 \bqn{ 
 	\lefteqn{	- \ln \pk{ \vk X(t_i) \le \vk x_i, 1 \le i\le n}}\notag\\ 
 		 	&=&	  \sum_{1\le l \le n} \E*{   \norm{ \vk Z(0)/\vk x_{l}}^\alpha _\IF
 			\mathbb{I}( infargmax_{  1 \le j \le n }   \norm{ \vk Z(t_j-t_l)/\vk x_{j}}^\alpha_\IF =l  )  } \notag \\
 		 	&=&	  \sum_{1\le l \le n} \E*{   \kappa(\vk Z) \frac{\norm{ \vk Z(0)/\vk x_{l}}^\alpha _\IF}{\kappa(\vk Z)}
	\mathbb{I}\bigl( infargmax_{  1 \le j \le n }   \norm{ \vk Z(t_j-t_l)/\vk x_{j}}^\alpha_\IF/\kappa(\vk Z) =l  \bigr)  } \notag \\ 
	&=&	  \sum_{1\le l \le n} \E*{    {\norm{ \vk \Theta(0)/\vk x_{l}}^\alpha _\IF} \label{infargm}
		\mathbb{I}\bigl( infargmax_{  1 \le j \le n }   \norm{ \vk \Theta(t_j-t_l)/\vk x_{j}}^\alpha_\IF  =l  
		\bigr)  }.
 }
Note that $\vk X$ is regularly varying, when $\vk X$ has \cadlag\ sample paths and hence $\vk \Theta$ can be derived also as a weak limit (under conditions on $\kappa$), see \cite{PH2020,MartinE}.  Since in view of \cite{Htilt, KumeE} $\vk X$ is stationary iff \eqref{tcfN00}  holds for all $F\in \mathcal{H}_0$, see also \cite{Hrovje,klem}  the above derivation  implies the following  result: 
\BEL If $\kK$ is such that $\vk Z$ has non-negative components, then 
  $\kK$ is    shift-invariant iff   $\vk X$ is stationary. In further  \eqref{re} holds, then $\kK$ is shift-invariant iff the finite dimensional distributions (fidi's) of $\vk X$ are given by \eqref{infargm}.
  \label{propikuandej}
\EEL

In view of \cite{dom2016, Hrovje} $\pk{  \SZ< \IF}=1$ is related to the existence of so-called  Rosi\'nski   (or mixed moving maxima) representation of max-stable rf's, see also \cite{debicki2017approximation, PH2020}. 
Such representations are also referred to as mixed moving maxima (M3) representation \cite{kab2009, Hrovje} being closely related to dissipative/conservative decompositions, \cite{RZ, Genna04, Roy1,Roy, WangStoev, dom2016}.   The findings of \Cref{kabA} give several other equivalent conditions for such representations, which  are considered in detail in \cite{EHPara}. \\

  \subsection{Maximal indices}
Let $\EIE$ be  a  full rank lattice on $\TT$  with $\abs{\EIE}=\IF$ or $\EIE=\R^l=\TT$. Consider further a given  $\kk$   and define the maximal index 
\bqn{ \label{nu}
	\bobo_{\vk Z}^{\EIE}= 
	\limit{n} n^{-l}\E*{ \max_{t\in [0,n]^l \cap \EIE} \kappa(B^{-t}  \vk Z)}.
}	
We drop the superscript $\EIE$ when  it is equal to $\TT$. It is clear that by \eqref{tcfN}   
$$
\bobo_{\vk Z}^{\EIE}=\bobo_{\ZT}^{\EIE}, \quad \forall \ZT \in \kk.
$$
If  
the stationary max-stable rf $ X$ has representer $\kappa(B^{-t} \vk Z)^{1/\alpha}, t\in \TT$  and $\alpha$-\FRE marginals, then  
$\bobo_{\vk Z}^{\EIE}$ is closely  related to the  distribution of supremum of $ X$, namely  
utilising \eqref{kode} and \eqref{conditionC1}
\bqn{\label{LY}
	-   \ln  \pk*{ \sup_{t \in \cap [0,n]^l \cap  \EIE  } X(t) \le  r n^{l/\alpha} }
 =  \frac{1}{r^\alpha n^l}   \E*{  \sup_{t \in \cap [0,n]^l \cap \EIE } \kappa(B^{-t} \vk Z) }    \to  \frac 1{ r^\alpha}  	\bobo_{\vk Z}^{\EIE}
}
as $n\to \IF$ for all $r>0$. Hence in this case $\bobo_{\vk Z}^{\EIE} $ is the extremal index of $X$ with respect to $\EI$. 
 
 Calculation of the extremal index of a stationary regularly varying rf continues to be a topic of great research interest. Its existence in general case cannot be guaranteed and therefore $\bobo_{\vk Z}^{\EIE} $ is commonly referred to as the candidate extremal index, see \cite{Hrovje, kulik:soulier:2020}. Given its relation to Pickands constants and related applications in statistics, the derivation of different representations of the maximal indices is of particular interest, see \cite{BojanS, KW,stoev2010max,HBernulli, MR3745388, Hrovje,klem, BP,ZKE,Planic}.

\BS Let $\tau\inr $ be such that   \Cref{CR2} holds. If   $\EIE$ is  a full rank lattice on $\TTT$ and $\kk$ has local rf $\vk \Theta$, then 
\bqn{ \label{horse} 
  \pk{\mathcal{B}_{\EIE,\tau}(\vk Y)=\IF}=1 \iff \bobo_{\vk Z}^{\EIE}=0 \iff \bobo_{\vk Z}^{\TTT}=0  \iff 	\pk{\mathcal{B}_{\TTT,\tau}(\vk Y)=\IF}=1.
}
\label{interesante}
\ES

\def\EI{\EIE}

 \BRM \label{rmE}
\begin{enumerate} [(i)]
	\item 
 Note that 	$\bobo_{\vk Z}^{\TTT}=0 \iff \pk{\SY =\IF}=1$ follows from  \cite{Genna04,Genna04c,MR2384479,MR2453345} or directly from \eqref{am12}. 
 \item \cER{ In view of \Cref{kabA}, several other equivalent conditions can be added to \eqref{horse}}. 
	  \item   If  $\vk Z$ has a.s.\ sample paths in $\Dset$  and $\kappa(f)= \norm{f(0)}^\alpha$, then  the last two equivalences in \eqref{horse} for $\tau=0, l=1$ follow from 
	\cite[Lem 2.5, Thm  2.9]{PH2020}  and for $\EIE=\Z^l$  by \cite[Thm 3.8]{HBernulli}, see also \cite{kulik:soulier:2020} for the case $\EIE=\Z=\TTT$.
\item \label{newPF}  Some novel representations for  $\bobo_{\vk Z}^{\EIE}   $ 
can be derived    by combining  \eqref{peshk1}, \eqref{peshk2} and \eqref{pil}.
\end{enumerate}
\ERM

\def\EI{\EIE}	

 \def\ve{\varepsilon}

\def\BAA{\widetilde{\mathfrak{A}}}

\def\BAA{\widetilde{ \mathfrak{A}_\diad}}
\def\BBAA{\widetilde{ \AA}}
\def\JJ{I_\diad}
\def\JJk{I_\diad}

\section{Proofs}
\label{sec:proofs}
\prooflem{OM1} Let 
$p_N (t)>0,\forall t\in \diad$ be the pdf of some $\diad$-valued rv and set
\bqn{\label{eq:JJk} 
	\JJk(f)= \int_\diad  \kappa(B^{-t}f)p_N(t)\lambda(dt). 
}
By \eqref{ZY} and the fact that $p_N$ is a pmf we have a.s.\ 
\bqn{ 
	\label{eq:JJIF} 
\JJ(\ZT) \in (0,\IF), \quad \forall \ZT\in \kK.
} 
Applying  
 \eqref{boll} for all $G \in \mathcal{H}_\alpha$ 
yields 
 \bqny{ \E*{G(\vk Z) } &=&
 	\E*{ G(\vk Z) \frac{\JJ(\vk Z)}{\JJ(\vk Z)}  }
 	\\ 
 &=&
 	\sum_{t \in \diad} p_N(t)\E*{ \kappa(B^{-t} \vk Z) \frac{ G(\vk Z) }{\JJ(\vk Z)}  }
 	=:\sum_{t \in \diad} p_N(t)\E*{ \kappa(B^{-t} \vk Z) F(\vk Z) }
 	 	\\ 
 	&=&\E*{G(\ZT) },
 } 
  where we used that $F$ is $0$-homogeneous for the derivation of  last line above.  By assumption \eqref{jet} we obtain thus   \eqref{parrap:E}.\\
   \def\JJ{I_{\diad}^{\IF}}
   \def\JJk{I_{\diad}^{\IF}}
 Suppose next that $\eqref{parrap:E}$ is satisfied and set 
 $$\JJ(f)=\int_\diad \kappa_\IF(B^{-t}f)p_N(t)\lambda(dt).$$
  Again by  \eqref{ZY} a.s.\ \eqref{eq:JJIF} holds for   $\JJ(\ZT)$. As above 
  for all $F \in \mathcal{H}_0$ we obtain 
 (set $\Gamma(\vk Z_{\star})(t):=\vk Z^{+}(t)-\vk Z^{-}(t)=\vk Z(t)$ and recall $\vk x_\star=(\vk x^+, \vk x^-), \vk x\inr^d$)
 \bqny{ \E{ \kappa(B^{-h}\vk Z )  F(\vk Z)}&=&\E*{ \kappa(B^{-h}\vk Z)  \frac{\JJk(\vk Z_{\star})}{\JJk(\vk Z_{\star})}F(\vk Z)}
 	=
 	\sum_{t \in \diad} p_N(t)
 	\E*{  \kappa_\IF(B^{-t}\vk Z_{\star})   \frac{ \kappa(B^{-t} \Gamma(\vk Z_{\star}  )) }
 		{\JJk(\vk Z_{\star})} F(\Gamma(\vk Z_{\star}))}\\
 		&=:&
 	\sum_{t \in \diad}p_N(t)\E*{ \kappa_\IF( B^{-t}\vk Z_{\star})  G(\vk Z_{\star}) }
 	\ttt \sum_{t \in \diad}p_N(t)\E*{ \kappa_\IF( B^{-t} \ZT_{\star})  G(\ZT_{\star}) }\ttt
 	\E*{ \kappa( B^{-h}\ZT) F(\ZT)}, \quad \forall h \in \TT,
 }
where we used that $G\in \mathcal{H}_0$.  Now we justify the second last equality, i.e., we need to show that 
\bqn{ \label{anmeldung} 
	\E*{ \kappa_\IF( B^{-t}\vk Z_{\star})  G(\vk Z_{\star}) } = \E*{ \kappa_\IF( B^{-t}\ZT_{\star})  G(\ZT_{\star}) }, \quad 
	\forall t\in \TT, \forall G \in \mathcal{H}_0.
} 
It follows from \eqref{kode} that the max-stable rf $\vk X_{\star}$ with representer $\vk Z_{\star}$ has also representer $\ZT_{\star}$. Since the law of $\vk X_{\star}$ is determined by the tail measure $\nu_{\vk Z_{\star}}$ (recall definition \eqref{nuZ}), which is  equal with $\nu_{\ZT_{\star}}$, then 
\cite[Prop 3.6, Rem 3.11,(ii)]{MartinE} implies that \eqref{boll} is satisfied with  $\kappa(B^{-h}\vk Z),\kappa(B^{-h}\ZT)$ substituted by $\kappa_\IF(B^{-h} \vk Z_{\star}),\kappa_\IF(B^{-h} \ZT_{\star}),$ respectively and thus \eqref{anmeldung} follows  establishing the proof.  
\QED

\def\JJ{I_\diad}
\def\JJk{I_\diad}
\def\IIj{I_\diad}

\prooftheo{lemAHDOOB2} 
		\eqref{boll}$ \implies $ \eqref{boll3}:
		The claim is clearly valid when $\TT=\Z^l$. Consider therefore the  case $ \TT=\R^l$. Let   $g \in \mathfrak{G}$, i.e.,  $g:\TTT\mapsto [0,\IF)$ is locally bounded and  $\lambda$-measurable and write
		\begin{equation}\label{rkn}
			R_{c,n}(f,g)= \frac{(2c)^l}{n^{l}}\sum_{ t\in (\Z/n)^l \cap [-c,c]^l}  \kappa(B^{-t} f)^\xi g(t), \quad f\in \Dset,
		\end{equation} 
		with $\xi \ge 0$.   Set  next   $\II  (f, g) = \int_\TTT \kappa(B^{-t} f)^\xi  g(t)\lambda(dt)$  and define 
		$$ F(f)=    \frac{ \Gamma_\xi  ( f )}{ \II  (f, g)   }, \quad  \Gamma_\xi\in   \mathcal{H}_\xi, f\in \Dset. $$
 \eB{We consider only this case for $F$, the other cases follow with the same arguments. }
Note in passing that  $F(\ZT), \ZT \in \kK$ is a well-defined rv (recall $\kappa(B^t \vk Z), t\in \TT$ is jointly measurable). By \eqref{ZY}, applying \Cref{lemGK} taking $d=1,\mathcal{A}=\R^l, \vk U(t)=\kappa(B^{-t}\ZT)^\xi, \gamma=z=0$ and $g_1(s)=s^\alpha, s\ge 0, g_2(t)=1, t\in \TT$ therein, we obtain 
\bqn{ \label{pimdot}
	\pk{\II(\vk Z,  g)>0}=	\pk{\II( \ZT,  g)>0}=1, \quad \forall \ZT \in \kK.
}
In fact, even when $\pk{\II(\vk Z,  g)=0}>0$, then by \Cref{lem:Belalp} 
$\pk{\II(\ZT,g)=0}>0$ and thus \eqref{boll3} obviously holds since both sides are equal to infinity. \\
We can assume without loss of generality that  $F$ is bounded since $F(f) \ind{F(f) \le n}, \eB{n\in \mathbb{N}}$ is also $0$-homogeneous and then we can apply the dominated convergence theorem for unbounded case. 
Hence  by \eqref{pimdot} the rv's $F(\vk Z),F(\ZT)$ are  strictly positive and finite. Clearly, the following map  
$$\mathfrak{U}_{k,n}(f)= \Gamma_\xi (f) /R_{k,n}(f,g), \quad  f\in \Dset , \ k,n> 0$$ 
 is a $\AA/\borel{\R}$-measurable belonging to  $\mathcal{H}_0$ for all large integers $k,n$ (we set $\mathfrak{U}_{k,n}(f)=0$ if $f$ is not integrable and interpret $0/0$ as $0$).  
For  all $h\in \TT,r>0$ and all positive integers $k,n$ large enough  by \eqref{tcfN}  for all $ \WZ\in \kK$
\bqn{\label{radio} 
	 \E*{ \kappa(B^{-h} \vk Z)  \mathfrak{U}_{k,n}(\vk Z)\ind{  \mathfrak{U}_{k,n}(\vk Z) \le  r} } &=&   
	\E*{ \kappa(B^{-h} \ZT)  \mathfrak{U}_{k,n}(\widetilde{ \vk Z})\ind{  \mathfrak{U}_{k,n}(\widetilde{ \vk Z}) \le  r}}<\IF.
}
In view of \eqref{conditionC1},    \Cref{lemAHDOOB}  yields for some sequence $c_k, n_{k}, k\inn$ 
$$ 
R_{c_k,n_{k}}(\ZT, g) \toas \II( \ZT, g), \quad \forall \ZT \in \kK$$
as $k\to \IF$. Hence for almost all  $r>0, \ZT \in \kK$ we have the convergence 
in distribution
 $$\mathfrak{U}_{c_k,n_{k}}(\widetilde{ \vk Z})\ind{ \mathfrak{U}_{c_k,n_{k}}(\widetilde{ \vk Z}) \le r} \todis 
\mathfrak{U}(\widetilde{ \vk Z})\ind{  \mathfrak{U}(\widetilde{ \vk Z}) \le  r}, \quad 
k\to \IF,
$$
where   
$\mathfrak{U}(\ZT) = \Gamma_\xi(\ZT)/\II( {\ZT}, g).$ 
Consequently, by \eqref{radio} 
\bqn{  \E*{  \kappa(B^{-h} \vk Z)   \mathfrak{U}(\vk Z)\ind{ \mathfrak{U}(\vk Z) \le  r} } &=&   
	\E*{  \kappa(B^{-h} \ZT)  \mathfrak{U}(\widetilde{ \vk Z})\ind{  \mathfrak{U}(\widetilde{ \vk Z}) \le  r}}< \IF.
\label{oriz}
}	

In view of   \eqref{pimdot} a.s.\ $\mathfrak{U}({ \vk Z}) ,\mathfrak{U}({ \ZT})\in [0,\IF)$. This is important since we can now apply monotone convergence theorem letting $r\to \IF$ to obtain 
\bqny{  \E*{ \kappa(B^{-h} \vk Z)    \mathfrak{U}(\vk Z)  } &=&   
	\E*{ \kappa(B^{-h} \ZT)    \mathfrak{U}(\widetilde{ \vk Z}) }\le \IF
}	 
establishing \eqref{boll3}.  \\ 
\eqref{boll3}$ \implies $ \eqref{boll1}: With the notation and the arguments of the proof in \Cref{OM1} for 
 $\ZT \in \kKY$  
applying the Fubini-Tonelli theorem and \eqref{boll3} for all $F \in \Hh_{\alpha }$ we obtain (recall \eqref{eq:JJIF})
 \bqny{\E{F(\vk Z) } &=& 
 	\E*{   \frac{ \JJ(\vk Z)}{ \JJ(\vk Z)}   F(\vk Z) }
 	=  \sum_{h\in \diad}   p_N(h)\E{ \kappa( B^{-h} \vk Z)    F(\vk Z)/\JJ(\vk Z)}  
 	 = \E{F(\ZT) }
 }
establishing \eqref{boll1}.  \\ 
If \eqref{boll1} holds, then clearly \eqref{boll} is satisfied,  hence the proof is complete. \QED

\prooflem{shpejt} $\kK$ is shift-invariant $\iff$ \ref{item:shpejt:0}, \ref{item:shpejt:00}: If $\kK$ is shift-invariant,  then 
\Cref{item:shpejt:0}-\Cref{item:shpejt:00} follows straightforwardly. Suppose next  that $\mathscr{C}_\kappa[\vk Z_\diad]$ is shift-invariant. Note that \eqref{conditionC1} holds with $\TT$ substituted by $\diad$. Since $\diad$ is a separant for $\kappa(B^t\vk Z), t\in \TT$, then 
\eqref{conditionC1} follows. The latter together with the stochastic continuity and the non-negativity of $\kappa$ imply that for all $t\in \TT$ and 
some $c$ sufficiently large 
$$ \IF > \E*{ \sup_{s\in [- c,c]^l} \kappa(B^{-s} \vk Z)}\ge 
 \limit{n} \E*{\kappa(B^{-t_n} \vk Z)} = \E*{\kappa(B^{-t} \vk Z)},
 $$
 with $t_n\in \diad, n\ge 1 $  such that $\limit{n} t_n = t$ (recall $\diad$ is a dense countable set of $\TT$). 
  In order to finish the proof, we need to show that for   all $ F\in \mathcal{H},h\in \TT$
 
 \bqny{  \E*{ F( \vkT_h)} =\E{\kappa(B^{-h} \vk Z) F(\vk Z/ \kappa(B^{-h 	} \vk Z)^{1/\alpha} ) } 
 	&=&   
 	\E{\kappa(\vk Z) F(B^h\vk Z/ \kappa( \vk Z)^{1/\alpha} )}=\E*{ F(B^h \vkT) }  .
 }
  It suffices to show therefore that $\vkT_h$ and $B^h\vkT$ have the same fidi's. For $h\in \diad$ this is an immediate consequence of
 $\mathscr{C}_\kappa[\vk Z_\diad]$ being shift-invariant. For a general $h\in \TT$ it suffices to show that for all $t_1 \ldot t_n $ in $\TT$ and for all bounded continuous $F: (\R^d)^n \mapsto [0,\IF)$
 \bqn{\label{exam21}
 	\lefteqn{	  \E*{ \kappa(B^{-h} \vk Z) F( \vk Z(t_1) /\kappa(B^{-h} \vk Z)^{1/\alpha}  \ldot \vk Z(t_n) / 
 			\kappa(B^{-h} \vk Z)^{1/\alpha} )} }\notag \\
 	&=&   
 	\E*{ \kappa( \vk Z) F(B^h \vk Z(t_1)/\kappa(\vk Z)^{1/\alpha} \ldot B^h \vk Z(t_n) /\kappa(\vk Z)^{1/\alpha} ) }.
 }
  Since $\diad$ is dense in $\TTT$ there exists  $t_{ki},h_k\in \diad, k\inn$ such that 
 $$\limit{k} t_{ki}=t_i\in \TTT, \ 1 \le i\le n, \ \limit{k} h_{k}=h.$$
By  the stostochastic continuity of $\kappa(B^t \ZT), t\in \TT$ and $t_{ki}- h\in \diad$ (since $\diad$ is an additive group) using   the dominated convergence theorem (recall also \eqref{conditionC1}) the claim follows. \\
 \Cref{item:shpejt:00} $\implies$ $\kK$ is shift-invariant follows with similar arguments. 
  
\ref{item:shpejt:0}$ \implies $\Cref{item:shpejt:2}:  
From \eqref{eq:JJIF}  a.s.\   $I_\diad(B^N\vk Z)\in (0,\IF)$.
 Consequently, the following rf  
$$ \vk Z_N(t) =  \Bigl(\frac{\kappa(\vk Z)}{  I_\diad(B^N\vk Z) }\Bigr) ^{1/\alpha} B^N \vk Z (t), \quad t\in \TT$$
is well-defined and stochastically continuous.     By the assumption and the fact that $\diad$ is dense in $\TT$ it follows that 
\eqref{tcfN} holds for all  $h\in \diad, F\in \Hh_{\alpha}$. Next, given  $F\in \Hh_{\alpha}$,  by the  independence  of $N$ and $\vk Z$,  
 using further  \eqref{tcfN}  for the derivation of the third last equality below  we obtain  
\bqny{ 
\E{ F(\vk Z_N) } &=& 
\E*{  \intDD 
		 	\frac{ \kappa( \vk Z)}
		{ \JJ(B^h \vk Z)  }  F(     B^h \vk Z)  p_N(h) \lambda(dh)}  \\
	&=:&  \intDD \E{     \Gamma( B^h \vk Z) } p_N(h) \lambda(dh)\\
	&=&  \intDD \E{    \Gamma(\vk  Z) } p_N(h) \lambda(dh)\\
	&=&  \E*{    F(\vk Z)   \intDD   \frac{ \kappa( B^{-h}\vk Z)   p_N(h) } { \JJ(\vk Z)  } \lambda(dh)}  \\
	&=&  \E{ F(\vk Z)},
} 		 
hence $\vk Z_N \in \kK$.\\
\Cref{item:shpejt:0}$ \implies $\Cref{item:shpejt:4}: The proof is the same as above substituting $\diad$ by $\TT$ and the counting measure on $\diad$ by the Lebesgue measure on $\R^l$ if $\TT=\R^l$.\\
\ref{item:shpejt:2} or \ref{item:shpejt:4}   $\implies$  \ref{item:shpejt:0}: Clearly, if $\vk Z_N \in \kK$ for all   $\diad$-valued rv $N$ with  $p_N \in \pclas$, then it follows easily that $B^h\vk Z_N  \in \kK$ for all $h\in \diad$ implying that $\mathscr{C}_\kappa[\vk Z_\diad]$ is shift-invariant. The proof is the same when $p_N \in \pclasT$.\QED

\prooflem{theoA}   
Since by the assumption $\vk Z_N \in \kK$, by \Cref{shpejt}, \Cref{item:shpejt:0}  the claim follows  if we show that  
\bqn{ \label{pinocio} 
	\E{  \kappa(B^{-h}  \vk Z_N) F( \vk Z_N )}  
	=\E*{ F(B^h \vk \Theta) },  \quad \forall  	F \in \Hh_0, h\in \diad.
}
If $\vk \Theta$ is a spectral tail rf, then in view of \eqref{eqDo20}	and the assumption that $\pk{ \kappa(\vk \Theta)=1}=1$ by    the Fubini-Tonelli  theorem for all $F \in \Hh_0, h\in \diad$ we obtain
 	\bqny{ 
		\E{  \kappa( B^{-h} \vk Z_N) F( \vk Z_N )}  
		&=&  \intDD \E*{  \kappa(B^{y-h}\vk \Theta) \kappa(\vk \Theta)\frac {  F( B^{y}\vk \Theta )}
			{ \JJ( B^{y}\vk \Theta ) }  } p_N(y) \lambda(dy) \notag \\ 
		&=& \E*{ \frac {F(B^h \vk \Theta )}	{\JJ( B^{h}\vk \Theta ) }  \intDD    \kappa( B^{h-y}\vk \Theta)   			p_N(y) \lambda(dy)}\notag\\
		  &=&  \E{ F(B^h \vk \Theta )},
	}	
	 hence \eqref{pinocio} holds. When $p_N \in \pclasT$ the proof follows with the same arguments as above.\\
	The converse is consequence of \eqref{tcfN} and \Cref{shpejt}. 
	\QED

\prooftheo{altFiorentinischeTrad} 
From \Cref{remalso}, \Cref{tYYY} it follows that  $\kappa(\vkT)^{1/\alpha}=R$ is an 
$\alpha$-Pareto rv and   $\vk \Theta=\vk Y/R$ is independent of $R$. 
\cER{Since $\vk Y$ satisfies \eqref{justen}, then $\vk \Theta$ also satisfies \eqref{justen}.} 
Next, by \Cref{remalso}, \Cref{tYYY} and the $\alpha$-homogeneity of $\kappa$     (recall $\lambda_\alpha(dz)=\alpha z^{-\alpha -1} dz $))
\bqny{ \lefteqn{ \int_x^\IF \E{ \Gamma(  z\vk \Theta) 
		\ind{ \kappa( z B^{-h}\vkT )> 1}} \lambda_\alpha(dz)}\\
	  &=&x^{-\alpha}\int_1^\IF \E{ \Gamma(   xz\vk  \Theta ) 
		\ind{ \kappa(  x zB^{-h}\vkT )> 1}} \lambda_\alpha(dz) \\
&=&
	x^{-\alpha}\E{  \Gamma (  x\vk Y) \ind{ \kappa( x B^{-h} \vk Y) >1  }     }\\
&=&\E{  \Gamma (  B^h \vk Y) \ind{ \kappa(  B^{h} \vk Y) >x^\alpha  }     }
, \quad 
	\forall \Gamma \in \wwH, \forall h\in \TT, \forall x>0.
}  
Consequently, letting $x\downarrow 0$  the monotone convergence theorem yields \eqref{tYYGen}.  
For all $  F\in \Hh_0,h\in \TTT $ by \Cref{remalso}, \Cref{tYYY} and the Fubini-Tonelli theorem using further $\vk Y= R \vk \Theta, \pk{\kappa(\vk \Theta)=1}=1$ and  the $\alpha$-homogeneity of $\kappa$ we obtain 
	\bqny{\E{ \kappa( B^{-h}\vk \Theta) F(\vk   \Theta  )    } 
		&=& \E*{  \kappa(  B^{-h} \vk Y )/\kappa(\vk Y)  
			 F(  \vk Y )} \notag \\
		&=& \E*{ z^\alpha\frac{F( z \vk Y )}{ \kappa(  z \vk Y) } \int_0^\IF  \ind{ \kappa( z  B^{-h} \vk Y)> 1}	\lambda_\alpha(dz)   } \notag \\
		&=&    \int_0^ \IF  z^{2\alpha}  
	\E*{  \frac{ F( B^h\vk   Y )} {\kappa( B^h \vk Y)}   \ind{  \kappa( B^h \vk   Y/z)  >  1  } }  \lambda_\alpha(dz) \\
		&=&   \E*{ \frac{ F(  B^h \vk  Y)}{ \kappa( B^h \vk Y) }   \int_0^{ \kappa(  B^h \vk Y)^{1/\alpha} }    \drzz }  \\
		&=&  
		\E*{   { F(  B^h \vk  \Theta)}  \ind{ \kappa( B^h\vk \Theta) \not=0}}.
	}
Hence  $\vk \Theta$ is a spectral tail rf. 
Next,  let $\kk=\mathscr{C}_\kappa(\vk \Theta)$ be as in \Cref{remKA} and write $\vk Y= R \vk \Theta$, which is a tail spectral rf.
In view of \Cref{CR2} and \Cref{lemGK}  we have 
\def\IK{J}
$$\IK (\vk Y)= \int_{K} \kappa( B^{-t} \vk \Theta)^{\tau} \ind{ \kappa(B^{-t} \vk Y) > 1} \lambda(dt)< \IF, \quad K=[-c,c]^l \cap \TT.$$  
Since $\kappa(B^t \vk Y), t\in \TT$ is jointly measurable, $\IK(\vk Y)$ is well-defined and further 
$\pk{\kappa( \vk Y)>1}=1$  together with \Cref{lemGK} (taking $d=1, \mathcal{A}=K,\vk U=\kappa(\vk Y), g_1(s)=s^\tau,s\ge 0,\gamma=z=1$ and $g_2$ equal 1 therein) imply  $\pk{\IK(\vk Y)>0 }=1.$ Moreover, applying again \Cref{lemGK} for all $s>0$ on the event $\{   \sup_{t\in K} \kappa(sB^N \vk Y) > 1\}$ we have that 
$\IK( s B^{N-t} \vk Y)> 0$ and thus on that event  $\IK( s B^{N-t} \vk Y)/\IK( sB^{N-t} \vk Y)=1$ for all $t\in K$.    
Consequently, borrowing the arguments of \cite{PH2020}, using additionally \Cref{shpejt} and the $\alpha$-homogeneity of $\kappa$ (with $N$ and $p_N$ as therein)  we obtain
\bqny{ \E*{ \sup_{t\in K}  \kappa(B^{-t}\vk Z) }
	&=& 
	\E*{  \frac{ \sup_{t\in K} \kappa( B^{N-t}\vk Y )}{\int_{\TT} \kappa(B^{N-t}\vk Y ) p_N(t) \lambda(dt)}}\\
	&=& \E*{ \int_0^\IF  \frac{ \ind{  \sup_{t\in K} \kappa( sB^{N-t}\vk Y) > 1} }{\int_{\TT} \kappa(B^{N-t}\vk Y  ) p_N(t) \lambda(dt)} \lambda_\alpha(ds) }
	\\
	&=& \int_{K} \int_0^\IF \E*{ \frac{\kappa(s B^{N-x} \vk Y)^\tau \ind{ \sup_{t\in K} \kappa(s B^{N-t}\vk Y) > 1, \kappa(s B^{N-x}\vk Y) > 1} }{\IK(s B^N \vk Y)\int_{\TT} \kappa(B^{N-t}\vk Y) p_N(t) \lambda(dt) {\tiny }} }\lambda_\alpha(ds)  \lambda(d x)\\
	&=& \intT \int_{K} \int_0^\IF \E*{ \frac{s^\alpha \kappa( s B^{h-x} \vk Y)^\tau \ind{ 
				\kappa(s B^{h-x}\vk Y) > 1} }{\IK(s B^h \vk Y)
			\int_{\TT} \kappa(s B^{h-t}\vk Y)  p_N(t) \lambda(dt) } }\lambda_\alpha(ds)  \lambda(d x) p_N(h)\lambda(dh)\\
	&=& \intT \int_{K}  \E*{  \int_0^\IF\frac{\kappa( \vk Y)^\tau  \ind{    \kappa( B^{x-h}\vk Y/s) > 1} }
		{\IK( B^{x} \vk Y)\int_{\TT} \kappa(B^{x-t}\vk Y) p_N(t) \lambda(dt) }  \alpha s^{\alpha -1} \lambda(ds)}  \lambda(d x) p_N(h) \lambda(dh)
	\\
	&=& \int_{K}  \E*{ \frac {\kappa(\vk Y )^\tau} { \IK( B^x \vk Y)} \intT \frac{  \kappa( B^{x-h}\vk Y )  }{\int_{\TT} \kappa( B^{x-t}\vk Y)  p_N(t) \lambda(dt)}    p_N(h) \lambda(dh)}  \lambda(d x) 
	\\
	&=& \int_{K}  \E*{ \frac {\kappa( \vk Y)^\tau} { \IK( B^{x} \vk Y)}  }  \lambda(d x), 
}
 where we used \Cref{remalso}, \Cref{tYYY} (this is crucial for the proof) to derive the last third line above. Consequently, we  have  
\bqn{ \E*{\sup_{t\in K} \cEE{\kappa(B^{-t}\vk Z) }}  &=& 
	\int_{  K}\E*{ \frac{ 1}{\int_{  K }  \kappa(B^{x-s}\vk \Theta)^\tau \ind{ \kappa( B^{x-s}\vk Y)> 1 }  \lambda(ds)}}
	\lambda(dx) \in (0,\IF)
	\label{nota}
}	
 implying that     $\vk Z \in \mathscr{C}_\kappa(\vk \Theta)$ satisfies \eqref{conditionC1} is equivalent with \eqref{ey}. \\ 
It is clear that \eqref{tYYGen} implies \eqref{eqDo20} for all    $\Gamma\in  \wHO,  h\in \TT$ and 
 we can write again $\vk Y= R \vk \Theta$ with $R$ independent of $\vk \Theta$ being $\alpha$-Pareto. Indeed, if $\vk Y$ is a spectral tail rf, then \eqref{nota} follows  as mentioned above. Consequently, we can define an $\alpha$-homogeneous rf generated by $\vk Z_N$ as in \Cref{shpejt}, which is shift-invariant. It follows that its local rf is $\vk \Theta$ establishing the proof.  \QED

\prooflem{domodosola}  By the assumption that $\vk Z$ satisfies \eqref{peres1}, then a.s.
$$ \{   \SZq < \IF\} \subset \{ \mj_1(\vk Z) \in \EI\} .$$
 Taking $F(\vk Z)= \ind{\mathcal{S}_{\EIE} (\vk Z) =\IF}$ and applying \eqref{tcfN} (recall 
  $\E{ A;B}$ stands for  $\E{A \ind{B}}$ and \Cref{A1})
 \bqn{	 
 	\E*{   \mathcal{S}_{\EIE} (\vk Z)  F(  \vk Z) ;    \mj_1(\vk Z ) =0 } 	&=& 
 	\sum_{s\in \EI}\E{    \kappa(  B^s \vk  Z) F( \vk Z) ;    \mj_1(\vk Z )  =0 } \notag\\
 	&=&   \sum_{s\in \EI} \E{   \kappa(   \vk  Z)   F(   \vk Z) \ind{ \mj_1 (B^{-s}\vk Z)  =0 } } \notag\\
 	&=& \E{   \kappa( \vk  Z) F(    \vk Z)  \sum_{s \in \EI } \ind{ \mj_1 (\vk Z)  =s} } \notag\\
 	&=& \E{\kappa( \vk Z) }\E{ F(  \vk \Theta );\mj_1(\vk \Theta) \in \EIE   }. \label{peshk1}
 }
 Since  $\E*{ \mathcal{S}_{\EIE} (\vk Z)F(\vk Z) ; \mj_1(\vk Z) = \EIE    }\in \{ 0, \IF\}$ we have 
$$ 
\E{ F(  \vk \Theta );\mj_1(\vk \Theta) \in \EIE   }=\pk{ \STq < \IF, \mj_1(\vk \Theta) \in \EIE   }=0.
$$
Hence   \Cref{lemAH} implies 
 $$\{ \mj_1(\vk Z) \in \EI\}  \subset \{  \SZq < \IF\}.$$
  Taking $\mj_1$ equal to  the infargsup map 
$\Ifasup$ establishes thus \eqref{ideT00}. As above, substituting $\mj_1$ by $\mj_2$ we have
$$0=\pk{ \STq < \IF, \mj_2(\vk \Theta) \in \EIE   }=\pk{ \mathcal{S}_\EI(\vk Y) < \IF, \mj_2(\vk Y) \in \EIE   }=0, $$
hence by   \eqref{peres2}  modulo null sets $ \{\mathcal{S}_\EI(\vk Y) < \IF\}=   \{\mj_2(\vk Y) \in \EIE\}  $ and thus the first two equalities in \eqref{ideT0} follow for $\mathcal{I}=\mj_2$.  Next since  $K \cap \EI$ has only finite number of elements for  every compact $K\subset \R^l$, then 
\bqn{ \label{sars} \quad \quad\mbox{}\{ \mathcal{S}_\EI(\vk Y)< \IF \}= \{  \STq < \IF  \}  \subset \Bigl \{ \lim_{ \sum_{i=1}^l \abs{t_i} \to \IF, t\in \EI} \kappa( B^{-t}\vkT) =0   \Bigr \} \subset
\{\JJqx < \IF \}
}
holds modulo null sets. 
Further, we have a.s.
$$ 
	\MME= \sup_{t\in \EI} \kappa(B^t \vk Y)^{1/\alpha} \ge \kappa(  \vk Y)^{1/\alpha} >1.
$$
Since by \Cref{lemKK} $ \JJqx< \IF  $ implies $ \MME < \IF $, then  
   in view of \eqref{fU} and \Cref{lemFUND}
$$ 
\{\JJqx < \IF\} \subset \{ \mathcal{S}_\EI(\vk Y) < \IF\} $$
 modulo null sets implying 
\bqn{ \label{then} 
	\{ \mathcal{S}_\EI(\vk Y)< \IF \}= 	\{\JJqx < \IF \}.
}
Let  $\Gamma \in \Hh_0$ be  shift-invariant with respect to $\EI$ and recall 
 $\JJqx= \sum_{i \in \TTd}  \kappa( B^i \vk Y)^{\tau}  \ind{\kappa( B^i \vk Y)> 1} $. Using  the Fubini-Tonelli theorem      
 for any $\Gamma \in \wwH$  and \Cref{A1},\Cref{A2} (recall \Cref{remalso}, \Cref{remalso:2})  
\bqn{ 
	\E*{   \JJqx \Gamma( \vk Y);  \mj_3(  \vk Y  ) =0 } 
	&=&  \sum_{i \in \TTd}   \E*{    \kappa(  B^i\vk Y) ^{\tau} \ind{\kappa( B^i \vk Y)> 1} \Gamma(B^i \vk Y); \mj_3( \vk Y  ) =0 } \notag \\
	&=&  \sum_{i \in \TTd}   \E*{    \kappa(  \vk Y)^{\tau} \ind{\kappa( B^{-i}\vk Y)> 1} \Gamma(\vk Y); \mj_3( B^{-i}\vk Y  ) =0 }\notag  \\
	&=&  \sum_{i \in \TTd}   \E*{    \kappa(  \vk Y) ^{\tau} \Gamma(\vk Y) ; \mj_3( \vk Y  ) =i }\notag \\
	&=&    \E*{   \Gamma( \vk Y) \kappa(  \vk Y)^{\tau} ;  \mj_3(\vk Y) \in \cER{\EI}  }. \label{peshk2}
}
Taking $\Gamma(\vk Y) = \ind{ \JJqx  =\IF}/R^{1/\tau}$ we obtain 
\cER{$$	\E*{   \JJqx \Gamma( \vk Y);  \mj_3(  \vk Y  ) =0 } =
\E*{   \kappa(  \vkT)^{\tau} \ind{ \JJqx  =\IF} ;  \mj_3(\vk Y) \in \cER{\EI}  } \le 1.
$$
Since $\pk{\kappa(\vkT) =1}=1$, then }
$$\pk*{ \JJqx = \IF,  \mj_3 (\vk Y)  \in \EIE }=0$$
 and thus by  
\eqref{peres2}
 $$ 
\{ \mj_3(\vk Y) \in \EI\} = \{ \mathcal{B}_{\EIE,0} (\vk Y )< \IF\},
$$
hence the proof follows utilising further \eqref{sars}. \QED 

\prooftheo{kabA}   
Applying \Cref{lemGK} establishes \eqref{am0}. 
Assume next that  $\pk{\mathcal{S}(\vk Z) < \IF}=1$, which in view of \Cref{lemAH} is equivalent with $\pk{\mathcal{S}(\ZT) < \IF}=1$.  Again by   \Cref{lemAH} 
$$q_1=\pk{ 0< \ST < \IF}=1.$$
Taking  $N$ independent of $\vk Z$ with pdf $p_N \in \pclasT$ we have 
\bqn{\label{eheq}
	\vk Z_N =  \Bigl( \frac{\kappa( \vk Z) }{  p_N(N)\mathcal{S}(\vk Z)  }\Bigr)^{1/\alpha} B^N \vk Z 
}
satisfies \eqref{tcfN}.   
Hence  it follows that   
$$
 \E*{ \sup_{t\in [-c,c]^l \cap \TTT} \kappa(B^{-t}\vk Z_N )} < \IF, \ \ \forall c \in (0,\IF),
$$  
which yields 
$$ q_2=\pk*{ \intT \sup_{t \in  [-c,c]^l \cap \TTT  } \kappa(B^{s-t}\vk \Theta) \lambda(ds) < \IF }=1, \quad \forall c\in (0,\IF)
$$
and thus by \Cref{lemFUND}    
\bqn{\label{bo1} 
	\{  \ST< \IF \} \subset  \Bigl\{ \intT \sup_{t \in  [-c,c]^l \cap \TTT  } \kappa(B^{s-t}\vk \Theta ) \lambda(ds) < \IF \Bigr\}, \quad \forall c\in (0,\IF).
}
In fact the above is an equality since the reverse inclusion clearly holds.  
When $q_2=1$, then also 
$$q_3= \pk*{ \lim_{ \sum_{i=1}^l \abs{t_i} \to \IF} \kappa(B^{-t}\vk \Theta ) =0 }=1.$$
Moreover if $q_3=1$, then   
$$
\pk{\mathcal{B}_{\EIE, \tau}(\vk Y) < \IF} =1.
$$  
Hence  \Cref{lemFUND} implies 
\bqn{\label{bo2} 
	 \Bigl\{ \intT \sup_{t \in  [-c,c]^l \cap \TTT  } \kappa(B^{s-t}\vk \Theta) \lambda(ds) < \IF \Bigr\} \subset \Bigl\{  \lim_{ \sum_{i=1}^l \abs{t_i} \to \IF} \kappa(B^{-t}\vk \Theta) =0\Bigr\} \subset  \{\mathcal{B}_{\EIE, \tau}(\vk Y) < \IF\}. 
}
In view of \eqref{horse},   \Cref{lemKK} and \Cref{lemFUND}  
 $$
 \{\BBYt < \IF\}= \{\mathcal{B}_{\EIE, \tau}(\vk Y) < \IF\}, \quad \{\BBYt < \IF\} \subset \{\ST  < \IF\},
 $$ 
 which together with \eqref{bo1}, \eqref{bo2} and \Cref{domodosola} establishes   
 \eqref{am10}-\eqref{am12}. \\ 
 The proof of \eqref{am20} and \eqref{am20c} follows with the same arguments using \Cref{lemAH} and  \Cref{easy}. \QED
	
\prooflem{propikuandej} If $\kK$ is    shift-invariant, then  the stationarity of $X$ follows from \eqref{tcfN0} and \cite[Thm 6.9]{Htilt} and the converse follows from the latter result. 
\QED  
	
\proofprop{interesante} First note that $ 		\bobo_{\vk Z}^{\EIE}   $  exists and is finite which follows from the \shiftinvariance of $\kk$ and the subadditivity of supremum, see also \cite{debicki2017approximation}. 
If    $\EIE$ is  a full rank lattice on $\TTT$ or $\EIE=\TTT$,  
	then by \eqref{nota}  
\bqn{ \label{pil}
		\bobo_{\vk Z}^{\EIE} =\frac{1}{\Delta(\EIE) }
		\E*{  \frac{ 1}{\intT \kappa(B^{-s}\vk \Theta)^\tau \ind{ \kappa( B^{-s}  \vk Y ) >1 } \lambda(ds)} }
		< \IF
	}	    
for all $\tau\inr$ such that \Cref{CR2} holds.
\def\EI{\EIE}
In our notation  $\Delta(\R^l)=1$ and if $\EIE$ is a full rank lattice on $\TTT$, then 
$\Delta(\EIE) =\abs{\det(A)}>0$ as in \eqref{bMat}, where  $A$ is  a non-singular base matrix of the full rank lattice $\EI$ i.e.,  $\EI=\{ A x, x\in \Z^l\}$. Consequently, $\bobo_{\vk Z}^{\EI}=0$ iff 
$$
	\pk*{ \JJqx = \IF}=1 
$$
and thus the first equivalence in \eqref{horse} is clear.\\
Suppose next that $\pk*{ \JJqx = \IF}=1$. By \Cref{domodosola}  
$\pk*{ \mathcal{B}_{\EIE,0}(\vk Y) = \IF}=1$. 
Since    for all $n\inn$ we have that $	\pk*{ \mathcal{B} _{ 2^{-n}\EIE, 0}(\vk Y ) = \IF}=1 $ and 
further  $2^{-n}\EIE$ is  also a full rank lattice on $\R^l$, then by the first equivalence in \eqref{horse}  and \eqref{3mace}  
we have that $\bobo_{\vk Z}^{\R^l}=0$, which is equivalent with $\pk{ \mathcal{B}_{\R^l,0 }(\vk Y)= \IF}=1$. 
If the latter holds, using that  
  $\bobo_{\vk Z}^{\EIE} \le \bobo_{\vk Z}^{\R^l}=0$, 
  then  
$\bobo_{\vk Z}^{\EIE}=0$ establishing \eqref{horse}.   
\QED

\def\JJxi{\mathbb{B}_\EIE(\vk Y/x)} 
\def\MPYn{ \mathcal{M}_{\EIE_n}(\vk Y)}
\def\JJqn{\mathcal{B}_{\EIE_n}(\vk Y)}
\def\TT{\EIE}
\def\TT{\mathcal{T}}

\section{Appendix}
Recall that in our notation convergence in probability and a.s.\ are denoted by $\toprob$ and $\toas$, respectively.
\label{sec:AAL}
\BT \label{lemGK} Let   $g_1:\R^d \mapsto [0,\IF]$ be  Borel measurable   and let    
$g_2:  \R^l \mapsto [0,\IF)$ be  Lebesgue measurable, almost everywhere positive and locally bounded. 
Let further    $\mathcal{A} \subset \R^l$  be open and  $\vk U\in \SCD$. If   
 $ \pk{\sup_{t\in \mathcal{A}} \norm{\vk U(t)}_\IF>z}>0, z\inr$, then for all $\gamma \in [-\IF, z]$
  \bqn{ \pk*{\sup_{t\in \mathcal{A}} \norm{\vk U(t)}_\IF>z, \int_{\mathcal{A}} g_1(\vk U(t)) \ind{ \norm{\vk U(t)}_\IF> \gamma  } g_2(t) \lambda(dt)=0 }=0,
  }
  provided that $g_1(f(t)) >0$ for all $t\in \mathcal{A}$ such that $\norm{f(t)}_\IF>\gamma, f\in \Dset	.$
\ET
\prooflem{lemGK} 	By the assumptions on $g_1$  we have that $g_1(\vk U(t)),t\in \R^l$ is measurable and non-negative. By the Fubini-Tonelli Theorem (see e.g., \cite[Thm 2.7]{Doob})
$\mathcal{I}_{\mathcal{A}}(\vk U)=\int_{\mathcal{A}} g_1({\vk U(t)} ) g_2(t)\lambda(dt)$ is a non-negative rv.\\
 For $d=1$ and $g_2$ constant the claim follows from \cite[Thm 2.1]{GeoSS}. The extension $d>1$ and $g_2$ non-constant follows with the same argument as the proof of the aforementioned result.
\QED 

\cER{
\BEL Let $X_{m,n},Y_{n}, n,m \inn$ be rv's defined on the same probability space. If 
$$X_{m,n} \toprob Y_m, \quad n\to \IF, \quad Y_m \toprob Y, \quad m\to \IF$$
holds, then there exists a non-decreasing sequence of integers $m_{k},n_k\inn$ such that $X_{m_k,n_k} \toas Y$ as $k\to \IF$.
\label{lem:phborsh} 
\EEL
}
\prooflem{lem:phborsh} Since the convergence in probability is metrizible
by the metric 
\bqn{\label{met} d(X,Y)= \E{\abs{X-Y}/(1+ \abs{X-Y}},
}
then by \cite[Lem A.1.3]{kulik:soulier:2020}
there exists a subsequence of non-decreasing integers $m_n, n\inn$ such that  $X_{m_{n_k},n_k} \toprob Y$ as $k\to \IF$ holds in probability. Hence there exists another a subsequence of non-decreasing integers $m_{k}, n_k\inn$
such that   $X_{m_{k},n_k} \toas Y$ as $k\to \IF$ establishing the proof.
\QED

\BEL Let   $g_1, g_2, \vk U$ be as in \Cref{lemGK}. If further 
\bqn{ \label{A} \int_{ \mathcal{A} }\E{ g_1(\vk U(t))}g_2(t) \lambda(dt)< \IF} 
  is valid for all $\mathcal{A}=[-c,c]^l, c>0$,   then        
	$\tilde R_{c_i,n_{i}}(g_1( {\vk U}),  g_2) \toas \mathcal{I}_{\R^l}(\vk U), i\to \IF$, where 
	$c_i,n_{i}, i\inn$	are two sequences of increasing positive integers converging to $\IF$ as $i\to \IF$ and 
	$$
	\tilde R_{c,n}(g_1(f),g_2)= \frac {\cEE{(2c)^l}} {n^l}\sum_{ t\in (\Z/n)^l \cap [-c,c]^l}  g_1(f(t))g_2(t), \quad f\in \Dset.
	$$
	\label{lemAHDOOB}
	\EEL 
	
	\prooflem{lemAHDOOB}
 Since  $\mathcal{I}_{[-c,c]^l }(\vk U) \toas \mathcal{I}_{\R^l}(\vk U), c\to \IF$ in view of \Cref{lem:phborsh} it suffices to show that 
	   $\tilde R_{c,n_k}(g_1(\vk U)g_2) \toas  \mathcal{I}_{[-c,c]^l}(\vk U), k\to \IF$   for a given $c>0$ and some non-decreasing sequence of integers $n_k,k\inn$ such that $\limit{k} n_k=\IF$. 
	The claim for $l=1$ is consequence of the derivations in \cite[p.\ 329-230]{MR636254}. Borrowing those arguments, leads to the proof for the case $l$ is a positive integer.
	\QED

	\BEL
Let $\kK$ be given and let $A \subset \R$  be a Borel set. 
If $G\in \mathcal{H}$ is such that $\ind{G(\cdot) \in A}\in \Hh_0$, then 
$$ \pk{G(\vk Z) \in A}>0  \text{ is equivalent with } 	\pk{G(\ZT) \in A}>0, \quad \forall \ZT \in \kK.$$
\label{lem:Belalp}	
\EEL	  
\prooflem{lem:Belalp} Using \eqref{eq:JJIF}, in view of 
\eqref{boll1}
\bqny{ 
	  \E{  \JJ(\vk Z) \ind{G(\vk Z) \in A }} &=&\sum_{t \in \diad} \E{  \norm{\vk Z(t)}^\alpha  \ind{G(\vk Z) \in A }} \\
	&  =&\sum_{t \in \diad} \E{  \norm{\ZT(t)}^\alpha  \ind{G(\ZT) \in A }}=
	\E{  \JJ(\ZT) \ind{G(\ZT) \in A }},
}
hence the claim follows. \QED

Recall  the Polish metric space $(\DDset,d_{\DDset})$ in \Cref{def:dc} and note that  the Borel $\sigma$-field of $\DDset $ agrees with the cylindrical $\sigma$-field $\AA$, see e.g., \cite{ferger2015arginf,MartinE}.\\
For $\vk X(t),t\in \TT$ with a.s.\ sample paths in $\DDset $ that is regularly varying with 
$\alpha$-homogeneous tail measure $\nu$, in view of \cite{kulik:soulier:2020,klem, MartinE} we have that $\nu=\nu_{\ZT}$ for some $\ZT \in \kK$, where  $\norm{\cdot}: \R^d \mapsto [0,\IF)$ is a norm on $\R^d$ and $\kappa(B^{\cdot}f)=\norm{f(\cdot)}^\alpha$.\\ 
Let $\mathcal{B}_0$ be a boundedness, which consists only of sets $A \in \AA$ such that  
for all $f\in A$ we have  $d_{\DDset }(f,0)>\ve_A$ for some $\ve_A>0$,  see \cite{kulik:soulier:2020,PH2020,MartinE}. 
Consider some positive sequence $a_n, n\ge 1$ such that $n \pk{ \vk X /a_n \in \cdot }$ converges weakly to $\nu_{\ZT}(\cdot)$ as $n\to \IF$ with respect to the boundedness $\mathcal{B}_0$, see \cite{kulik:soulier:2020,PH2020, MartinE} for more details.\\
 We formulate next a lemma on 1-homogeneous maps which is a minor extension of \cite[Prop 2.5]{jansen2022taildependence}, see also the related result \cite[Lem A.7]{stA}. In the following $H: \DDset  \to [0,\IF]$ is called lower semi-continuous if $A_x=\{ f \in \DDset: H(f)> x\}$
is an open subset of $\DDset $ for all $x>0$.  

\BEL Let $F\in \mathcal{H}_1$ \cED{with $F(0)=0$} be lower semi-continuous, non-negative and continuous at 0. If $\vk X$ is as above, then for all $x>0$
\bqn{ 
	\limit{n}  n \pk{ F(\vk X) > a_nx} = x^{-\alpha}\E{ F^\alpha(\vk Z)}< \IF.
}
\label{lemZJ}
\EEL 

\prooflem{lemZJ} 
 The continuity at zero of $H$ yields  that   $d_{\DDset}(0,f)> \ve$ for some $\epsilon >0$ and all 
 $f\in A_x,x>0$, which implies that for some hypercube $K\subset \R^l$ we have  
$$\sup_{t\in K \cap \TT} \norm{f(t)}  
> \ve$$
for all $f\in A_x$, see \cite{kulik:soulier:2020, MartinE}. Consequently, we have 
$A_x \in \mathcal{B}_0$. Since tail measures on $\AA$ are such that $\nu_{\ZT}(E)< \IF$ for all $E \in \mathcal{B}_0$, then by the $-\alpha$-homogeneity of $\nu_{\vk Z}$ and the 1-homogeneity of $F$
$$\nu_{\ZT}(A_x)= x^{-\alpha}\nu_{\ZT}(A_1)= x^{-\alpha}\E{ F^\alpha(\vk Z)} < \IF.
$$
Since $A_x$ is open its frontier is a subset of $A_x^{=}=\{ f \in \DDset: F(f)= x\}$, which also belongs to $\mathcal{B}_0$. By the $-\alpha$-homogeneity of $\nu_{\ZT}$ we have  that 
$$\nu_{\ZT}(A_x^{=})=0$$
 see also  \cite[Rem 3.2]{MartinE}) and hence the claim follows. See \cite[Prop 2.5]{jansen2022taildependence} or the derivation of \cite[(A.1)]{Stb}, where the last arguments appear for the case $d=1$ and $\TT$ has a finite number of elements, see also \cite[Lem 3.1]{Segers}.
\QED 

\BRM 
\begin{enumerate}[(i)]
		\item In view of \cite[Lem 2.2, Prop 2.1]{ferger2015arginf} 
		$$F(c,f)= \inf_{t\in [-c,c]^l} \norm{f(t)},  \quad  c\in (0,\IF]$$
		 satisfies the assumptions of \Cref{lemZJ};
	\item 
If $F_1,F_2$ and $X,Z$ are as in \Cref{lemZJ} such that $F_2(Z)>0$ a.s., then for $x,y$ positive 
\bqn{ 
	\limit{n}\pk{F_1(X)> a_nx \lvert F_2(X)> a_ny }
	=\frac{\E{ [\min(yF_1(Z), xF_2(Z)) ]^\alpha} }{\E{F_2^\alpha(Z)}}.
}
In particular, this is applicable for 
$F_1(c,f)= \sup_{t\in [-c,c]^l} \norm{f(t)}$ or $F_1(f)=\int_{[-c,c]^l} \norm{f(t)}\lambda( dt),$ with $c\in (0,\IF]$  and $F_2$ as in (i).
\end{enumerate}
\ERM

\BEL  
The local rf $\vk \Theta$ of $\kK$ is such that $\kappa(B^t \vk \Theta), t\in \TT$ is  stochstically continuous. 
\label{lemAS}
\EEL

\prooflem{lemAS} Let  $t_n \to t\in \TT$ as $n\to \IF$. 
By the $\alpha$-homogeneity of $\kappa$
	\bqny{ 
		C_n=\E*{ \frac{\abs{  \kappa( B^{t_n} \vk \Theta) - \kappa(B^t \vk \Theta) } }{1+
				\abs{  \kappa( B^{t_n} \vk \Theta) - \kappa(B^t \vk \Theta) } } }
		&=&
		\E*{\kappa( \vk Z)  \frac{\abs{  \kappa( B^{t_n} \vk Z) - \kappa(B^t \vk Z) }/\kappa( \vk Z) }  {1+\abs{  \kappa( B^{t_n} \vk Z) - \kappa(B^t \vk Z) }/\kappa( \vk Z)  }}\\
		&=&  
		\E*{\kappa( \vk Z)   \frac{\abs{  \kappa( B^{t_n} \vk Z) - \kappa(B^t \vk Z) }} { \kappa( \vk Z) +
				\abs{  \kappa( B^{t_n} \vk Z) - \kappa(B^t \vk Z) } }}. 
	}
	Since by assumption $\kappa(B^t \vk Z),t\in \TT$ is stochastically continuous, then  $\abs{  \kappa( B^{t_n} \vk Z) - \kappa(B^t \vk Z) } \toprob 0 $  as $n\to \IF$. Further, using that  $\kappa(\vk Z)$ is non-negative and $\E{ \kappa(\vk Z)}< \IF$ we obtain by the dominated convergence theorem $\limit{n} C_n=0$. 
	 Using that  the convergence in probability is metrizible (recall \eqref{met}), the claim follows. 
Note in passing that by  \cite[Thm 1, p.\ 171 \&  Thm 5, p.\ 169]{MR636254} both $\vk \Theta$ and $\kappa(B^t \vk \Theta), t\in \TT$ have  a jointly measurable   and separable version with separant $\diad$. 
\QED

\def\ko{K_{\vk 0}}

\BEL  
Let $\kk$ with local rf $\vk \Theta$ be given.  If   
  $F \in \wHO$ is non-negative, then  the following are equivalent: 
  \begin{enumerate}[(i)]
  	\item \label{elko:1}   \cEp{$\E{F(\ZT)} =0$ for some (and then all) $\ZT \in \kk$}; 
  	 	\item \label{elko:2} $\E{F(B^h\ZT)} =0$ for some (and then all) $\ZT \in \kk$ and for all  $h\in \TT$;
  	\item \label{elko:3} $\E{F(B^h\vk \Theta)} =0$ for all  $h\in \TT$.
  \end{enumerate}
Moreover, if  $F$ is  bounded by some constant $c>0$,   then  \Cref{elko:1}-\Cref{elko:3}  hold with $c- F$ instead of $F$.
\label{lemAH}
\EEL 

{
\prooflem{lemAH}  Recall the definition of  $\JJ(\vk Z)$ in \eqref{eq:JJk} and   let 
 $p_N \in \pclas $. \\
$\ref{elko:1} \implies \ref{elko:2}$: Since by the assumption $F \in \wHO$ is non-negative, then when \ref{elko:1} holds a.s.\ $F(\ZT)=0$. Using \eqref{tcfN} for all $h\in \TTT, \ZT^* \in \kk$   
	\bqny{  \E*{F( B^{h}\ZT^*  ) 	\JJk(\ZT^*)}=     
		\E*{F( B^{h}\ZT  ) \JJk(\ZT)}  =
	\intD{t}  \E*{  \kappa(B^{h-t} \ZT)F( \ZT  ) } p_N(t)  =  0.
}
Consequently, $\JJk(\ZT^*)\in(0, \IF)$ a.s.\ implies  $ F( B^{h}\ZT^*  ) =0$. 

$\ref{elko:2} \implies \ref{elko:3}$: For all $h\in \TT,\ZT \in \kk$   
$$\E{F(B^h\vk \Theta)} = \E*{  \kappa(\vk Z) F( B^h \vk Z/ \kappa(\vk Z)^{1/\alpha}) }= 
\E*{  \kappa(\vk Z) F( B^h \vk Z ) }=\E*{  \kappa(\ZT) F( B^h \ZT) }$$
 and thus 
\ref{elko:3} holds. \\
$\ref{elko:3} \implies \ref{elko:1}$: By 
the \shiftinvariance of \cEp{$\kk$}  
(using    \eqref{tcfN} to derive  the last  equality below) 
\bqny{ 0 	 
	&=&   \intD{t}  \E*{  F( B^{t}\vk \Theta  ) } p_N(t) = 
	 \intD{t}  \E*{  \kappa( \vk Z)  F( B^{t}\vk Z  ) } p_N(t)  =  \E*{  F( \vk Z  ) \intD{t}
	 	\kappa( B^{-t} \vk Z)  p_N(t)},
}
hence  \ref{elko:1} follows for $ \vk Z$.\\
Next, if $\E{F(\vk Z)} =c$ we have by the assumption that $\bar F=c- F \in \wHO$ is non-negative     and $\E{\bar F(\vk Z)}=0$, which by the above is equivalent with $\E{\bar F(\vk \Theta) }=0$, hence the proof is complete. 	
 \QED 
}

As an  application of \Cref{lemAH} we have the following   characterisation for such $\kk$'s.  
\BK Given $\kk$ such that $\E{\kappa(\vk Z)}=1$,  then the following are equivalent:
\begin{enumerate}[(i)]
	\item \label{elker:1} $\kappa( \vk Z)>0$ a.s.;
	\item \label{elker:12} For some (and then for all) $\ZT\in \kk$ we have  $\kappa(B^{h}\ZT)>0$ a.s.\  for all $h\in \TT$;
	\item \label{elker:2}  $\kappa(B^h\vkT)>0$ a.s.\ for all $h\in \TT$;
	\item \label{elker:22}  $ \vk \Theta\in \kk$.
\end{enumerate}
\label{lem:elker}
\EK
We state next  a fundamental property of shift-invariant maps.
\BEL \label{lemFUND}  Let $F_1, F_2  $ be two shift-invariant maps with respect to $\diad$ and 
\cER{let  $\setA \in \{\{0\},(0,\IF),\{\IF\}\} $.} 
\begin{enumerate} [(i)]
	\item 
	If $F_1,F_2  \in \wwH$ are such that for any tail  rf $\vk Y^* $   we have that  $\pk{ F_1(\vk Y^*) \in \setA}=1$ implies 
	$\pk{F_2(\vk Y^*) \in \setA}=1$, then 
	\bqn{\label{barfus}
		\{ F_1(\vk Y) \in \setA \} \subset  \{ F_2(\vk Y) \in \setA \}
	}	
	is valid for all tail rf's $\vk Y$; 
	\item 
	If $F_1,F_2 \in  \wHA$ are such that  for any spectral tail  rf $\vk \Theta^*$ we have that  $\pk{ F_1(\vk \Theta^*) \in \setA}=1$ implies 
	$\pk{F_2(\vk \Theta^*) \in \setA}=1$, then 
	\bqn{\label{barfus2}
		\{ F_1(\vk \Theta ) \in \setA \} \subset  \{ F_2(\vk \Theta ) \in \setA \}
	}	
	is valid for all spectral tail rf's $\vk \Theta $.
\end{enumerate}
\EEL

\BRM 
Utilising \Cref{lemAH} 
also the corresponding result of \Cref{lemFUND} for elements $\ZT \in \kk$ can be shown to hold. 
\label{easy}
\ERM	 

\prooflem{lemFUND}  Let $\vk Y$ be a tail rf with respect to $\kappa$ such that  
$\pk{F_1(\vk Y) \in \setA} \in (0,1]$ and define  ${\vk Y}^*= \vk Y \lvert F_1(\vk  Y)  \in \setA$, which by the assumption on $\setA$ and the shift-invariance of $F_1$ is a tail rf.
Clearly,  $\pk{F_1({\vk Y}^*)  \in \setA}=1$. By the assumption this implies that  $\pk{F_2({\vk  Y}^*) \in \setA}=1$ and hence 
$$1= \pk{F_2({\vk  Y}^*) \in \setA}= \pk{ F_2(\vk Y) \in \setA,  F_1(\vk Y) \in \setA}/
\pk{F_1(\vk Y)  \in \setA},
$$
which in turn yields   
$$\{ F_1(\vk Y) \in \setA \} \subset \{ F_2(\vk Y)  \in  \setA \}$$
modulo null sets establishing thus the first claim. The second claim follows with similar arguments.   
\QED

\def\intL{\int_{\EI}}
\BEL Let $\kk$ be given and let $\tau\inr$ be such that \Cref{CR2} holds. If $\pk{ \JJqx < \IF}> 0$ and $\EIE$ is a lattice on $\TT$  or  $\EIE= \TTT$, 
then  (recall $\MME= \sup_{t\in \EIE} \kappa( B^{-t } \vk Y)^{\cEE{1/\alpha}}$)
\bqn{ \pk{ 0<\JJqx <\IF , \MME=\IF}=0}
 and if further $\pk{\JJqx < \IF}=1$, we have 
\bqn{ \label{fU}
	\E*{\frac{\kappa(\vk Y)^\tau \mathcal{S}_\EI(\vk Y) }{ [\MME]^\alpha \JJqx }}=1.
}
\label{lemKK} 
\EEL 
\prooflem{lemKK} 
	Since by  \Cref{lemAS}  $\kappa(B^{-t} \vk Y)=R^\alpha\kappa(B^{-t} \vk \Theta),t\in \TT$ is stochastically continuous and hence belongs to $ \SCD$, then in view of \Cref{lemGK} (take therein $d=1, \mathcal{A}=\R^l, z=\gamma=1$, $\vk U= \kappa(\vk Y)^\tau$ and $g_1,g_2$ equal 1) 
	$$\MME=\sup_{t\in \EI} \kappa(B^t \vk Y)^{1/\alpha}> 1/z\  \implies \pk{ \mathcal{B}_{\EIE, \tau}(z\vk Y)> 0}=1, \quad z\in (0,1].$$
	  Moreover, $ \JJqx<\IF$ implies 
	$ \mathcal{B}_{\EIE, \tau}(z\vk Y)<\IF, z\in (0,1]$. Hence for all $z\in (0,1]$ by \Cref{remalso}, \Cref{tYYY} and the Fubini-Tonelli theorem (recall \Cref{remalso}, \Cref{remalso:2})  
	\bqny{  
		\lefteqn{ \E*{\frac{1}{\mathcal{B}_{\EIE, \tau}(\vk Y) }\ind{\mathcal{B}_{\EIE, \tau}(\vk Y) \in (0,  \IF), \MME = \IF   }}}\\
		&\le &\E*{\frac{ \mathcal{B}_{\EI, \tau}(z\vk Y) }{\JJqx  \mathcal{B}_{\EI, \tau}(z\vk Y)}
			 \ind{\JJqx\in (0, \IF), \MME >1/z   } }\\
		&= & 	\int_{\EI} \E*{\frac{\kappa(B^{-t}\vk Y)^\tau }{ \JJqx \mathcal{B}_{\EI , \tau}(z\vk Y) } \ind{\JJqx\in (0, \IF), z\MME > 1,   \kappa(z B^{-t}\vk Y) >1  }}\lambda(dt) \\
		&= & 	z^{\alpha}\int_{\EI } \E*{\frac{\kappa(z^{-1}\vk Y)^\tau }{ \mathcal{B}_{\EIE, \tau}(z^{-1}\vk Y) \JJqx } \ind{ \mathcal{B}_{\EI, \tau}(z^{-1}\vk Y)\in (0, \IF), \kappa(z^{-1}B^t\vk Y) >1  }}\lambda(dt)\\
		&= & 	z^{\alpha }\E*{\frac{\ind{\mathcal{B}_{\EIE, \tau}(z^{-1}\vk Y)\in (0,  \IF)} }{ \mathcal{B}_{\EIE, \tau}(z^{-1}\vk Y) \JJqx } 
			\int_{ \EI}  	\kappa(z^{-1}\vk Y)^\tau\ind{ \kappa(z^{-1} B^t\vk Y) >1  }\lambda(dt)}\\
		&= & 	z^{\alpha }\E*{\frac{\ind{\mathcal{B}_{\EIE, \tau}(z^{-1}\vk Y) \in (0,\IF)}\mathcal{B}_{\EI, \tau}(z^{-1}\vk Y)  }{ \mathcal{B}_{\cEp{\EIE, \tau}}(z^{-1}\vk Y)\JJqx}}\\
		&= & 	z^{\alpha}\E*{\frac{1  }{ \JJqx}}\\
		&<  & 	  \IF,
	}
	where the last inequality follows from  \eqref{ey}. Note that when $\EI$ is not equal to $\TT$ the above conclusion is obvious.  Consequently, letting $z\to 0$ yields
	$$\{\mathcal{B}_{\EIE, \tau}(\vk Y) < \IF \} \subset \{ \MME < \IF\}.$$
	For all $s\in (0,1)$ by \Cref{lemGK} we have 
	\bqny{ 
		\pk{ V_s>0, \MME < \IF}=1, \quad V_s=\mathcal{B}_{\EI , \tau}((s\MME )^{-1}\vk Y))
	}
	since  for all $s\in (0,1)$ on the event $\{ \MME < \IF\}$
	$$  
	\sup_{t\in \EI} \ind{ (s\MME )^{-1} \kappa( B^{-t}\vk Y)>1 }=1.
	$$
	By the Fubini-Tonelli theorem and (recall \Cref{remalso}, \Cref{remalso:2} and  $\lambda_\alpha(ds)=\alpha s^{-\alpha -1} ds $) 
	\bqny{  \int_0^1 
		\lefteqn{\E*{  V_s \ind{ V_s=\IF, \MME < \IF}  } 	\alpha s^{\alpha-1} ds}
		\\
		&=&   
		\int_{\EI}\int_0^1 \E*{  \ind{\kappa(B^{-h} \vk Y / (s\MME ))>1 }  \ind{ V_s =\IF, \MME < \IF}  } 	\lambda_\alpha(ds) 
		\lambda(dh)\\
		&=&   
		\int_{\EI}\int_0^\IF \E*{ \frac{1}{[\MME]^\alpha } \ind{\kappa(t^{-1}B^{-h} \vk Y)>1 }  \ind{ 
				\mathcal{B}_{\EI , \tau}(t^{-1}\vk Y)  =\IF, \MME < \IF}  } 	\lambda_\alpha(dt) 
		\lambda(dh)\\
		&=&   \int_{\EI}\int_0^\IF \E*{ \frac{1}{[\MME]^\alpha } \ind{\kappa(tB^h  \vk Y)>1 }  \ind{ 
				\mathcal{B}_{\EI , \tau}(\vk Y)  =\IF, \MME < \IF}  } 	\lambda_\alpha(dt) 
		\lambda(dh)\\
		&=&0
	}	
	under the  assumption  $\pk{\mathcal{B}_{\EI , \tau}(\vk Y)  =\IF}=0$, which we shall suppose next.
	Consequently, for all $s\in (0,1)$ up to a set with Lebesgue measure equal zero we have 
	\bqn{  \pk{ V_s\in (0,\IF) } =1.
		\label{aboOnline}
	}
	Using the latter implication, by the Fubini-Tonelli theorem  and (recall \Cref{remalso}, \Cref{remalso:2}, 
	$\IF \cdot 0$ as $0$) $\MME\in (0,\IF)$ a.s.\ since the origin belongs to $\EIE$ and we suppose that $\pk{\JJqx < \IF}=1$ and 
	$$ \mathcal{B}_{\EI , \tau}(z^{-1}\vk Y)=0 \implies \ind{\kappa(B^t\vk Y/z)>1}=0, \quad t\in \EI$$ 
	a.s.\ which are needed for justification of the third equality below (recall here \Cref{remalso}, \Cref{remalso:2})
	\bqny{ 
		\lefteqn{\E*{\frac{ \kappa( \vk Y)^\tau \mathcal{S}_\EI(\vk Y) }{ [\MME]^\alpha \JJqx }}}\\
		&= & \intL \E*{\frac{\kappa( \vk Y)^\tau \kappa( B^{-t}\vk Y)^\alpha  }{[\MME]^\alpha  \JJqx }}\lambda(dt) \\
		&= &\intL\int_0^\IF \E*{\frac{\kappa( \vk Y)^\tau  \ind{ \kappa(s B^{-t}  \vk Y)>1}  \ind{\JJqx\in (0,\IF)}} {[\MME]^\alpha   \JJqx }} \lambda_\alpha(ds)\lambda(dt)\\ 
		&= &\intL\int_0^\IF \E*{\frac{\kappa(B^t  \vk Y/\cEp{s})^\tau \ind{ \kappa( B^t \vk Y/s)>1} 
				\ind{ \mathcal{B}_{\EI , \tau}(s^{-1}\vk Y)\in (0,\IF)}  }{[\MME]^\alpha  \mathcal{B}_{\EI , \tau}(s^{-1}\vk Y)}} \lambda_\alpha(ds)\lambda(dt)\\ 
		&= &\intL\int_0^\IF \E*{\frac{\kappa(B^t \vk Y/(s\MME ))^\tau  \ind{ \kappa(B^t \vk Y/(s\MME) >1}
				\ind{ V_s\in (0,\IF)}  }{  V_s} } \lambda_\alpha(ds)\lambda(dt)\\ 
		&= &\int_0^1 \E*{\intL \kappa(B^t\vk Y /(s\MME ))^\tau   \ind{ \kappa( B^t \vk Y /(s\MME )) >1}  \lambda(dt)\frac{ \ind{ V_s\in (0,\IF)} }{ V_s } } \lambda_\alpha(ds)\\ 
		&= &1,  
	}
	where setting 1 in the upper bound of the integrand is justified by the fact that 
	$\ind{\kappa(B^{t} \vk Y /(s\MME))>1 }=0$ a.s.\ for all $s>1$ and further the $\alpha$-homogeneity of $\kappa$ was used for the derivation of the third equality.  
	\QED	  

\BEL If   $\TTT=\R^l$ and $\EIE$ is a full rank lattice  on $\R^l$, then  
\bqny{ \limit{n} 
	\limit{S} S^{-l}  
	\E*{   \Bigl[ \sup_{t\in [0,S]^l  } \kappa(B^{-t}\vk Z)  -   \sup_{t\in [0,S]^l \cap 2^{-n}  \EIE } \kappa(B^{-t}\vk Z)\Bigr]}=0.
}	 
Moreover,  with $\Delta(\EIE)$ defined in \eqref{bMat} we have
\bqn{  
	\limit{n} \frac{\bobo_{\vk Z}^{\EIE_n} }{\Delta(\EIE_n)}= \bobo_{\vk Z}^{\R^l} \in [0,\IF), 
	\quad \EIE_n= 2^{-n}\EIE. 
	\label{3mace}
}
\label{SSS}
\EEL
\proofprop{SSS} Assume for simplicity that $S$ is a positive integer and $\EI=\Z^l$. Set next $ L_n=2^{-n}\EI, n\inn $
	$$K_{\vk i}=\times_{k=1}^l [i_k,i_k+1)=[\vk i,\vk i+1), \vk i \in \Z^l, \ \  K_{\vk 0}=\times_{k=1}^l[0,1).$$
	For  all $n$ sufficiently large and any $\vk i \in \Z^l$ we have that $L_n\cap K_{\vk i}$ is non-empty. Note further that $L_n$ is a full rank lattice for any $n\inn$ and $\diad=\cup_{n\ge 1} L_n$ is a countable dense subset of $\TTT$ and thus it is a separant for  $\kappa(B^t \vk Z), t\in \TT$. Set below  
	$$ a(   K) = 
	\sup_{t\in K} \kappa(B^{-t}\vk Z) ,\ \ K\subset \TT.$$
	From  \eqref{conditionC1} and \eqref{tcfN} for all $n$ large we obtain 
	$$  
	\E*{   a( \ko \cap L_n) } =    \E*{ a(K_{\vk i} \cap L_n )}< \IF, \  \E*{   a( \ko ) }
	=    \E*{ a(K_{\vk i} )}< \IF, \  \forall \vk i \in \Z^l
	$$
	since   $\vk i \in L_n$ and $L_n$ is an additive group.
	Consequently,   for some positive integer $S$  
	\bqny{  
   \E*{   a(  [0,S)^l)  } -\E*{   a(  [0,S)^l \cap L_n) } 
		&=&  
		\E*{     \Bigl[  
			\max_{ \vk i \in [0,S-1]^l \cap \Z^l } a( K_{\vk i}  ) -	
			\max_{    \vk i \in [0,S-1]^l \cap \Z^l   } a(  K_{\vk i}\cap L_n )		 \Bigr]  } \\
		&\le &  
		\E*{     \max_{ \vk i \in [0,S-1]^l \cap \Z^l } \Bigl[  
			a( K_{\vk i}  ) -	 a(  K_{\vk i}\cap L_n )		 \Bigr]  } \\
		&\le &  	\sum_{\vk i \in [0,S-1]^l \cap \Z^l} 
		\E*{    \Bigl[  a( K_{\vk i})- a(K_{\vk i} \cap L_n ) 		\Bigr]  } \\
		&=& 	S^l\E*{  
			a( \ko) - a( \ko\cap   L_n ) }.
	}	
	Since $\kappa(B^{-t}\vk Z),t\in \TT$ is stochastically continuous  we obtain
\bqn{ \label{exampK} 
	M_n=\sup_{t\in \ko \cap L_n }	\kappa(B^{-t}\vk Z) \toprob \sup_{t\in \ko }   \kappa(B^{-t}\vk Z), \ \ 
	n\to \IF.
}	
Hence by the  dominated convergence theorem which is justified by \eqref{conditionC1} 
	\bqny{ \limit{n } \E*{    a( \ko) - a( \ko\cap   L_n ) 	 }
		=0
	}
	and thus the proof follows from \eqref{pil}.
	\QED

\BEL \label{lemGauss} If $\vk X(t)=(X_1(t)\ldot X_d(t)),t\in \TT$ is a centered $\R^d$-valued Gaussian rf  and let  $(Y, \vk X(t)), t\in \TT$ be jointly Gaussian defined on a probability space $(\Omega, \mathscr{F}, \mathbb{P})$.    
The law of $\vk X$ under the probability measure $\widehat{ \mathbb{P}}_v(A )= \E{ e^{ Y - v/2} \mathbb{I}(  A) }, A\in \mathscr{F}$ is the same as that of $\vk X + \vk C $ under $\mathbb{P}$, where 
$$\vk C(t)= (Cov( X_1(t), Y)
\ldot Cov( X_d(t), Y)), \ t\in \TT.
$$
\EEL 

	\prooflem{lemGauss}  
	Given  $t_1 \ldot t_n \in \TT$ we calculate  the df $\vk W=(\vk X  (t_1) \ldot \vk X(t_n))$ under  $\widehat{ \mathbb{P}}_v$
	by   exponential tilting of $(\vk X(t_1) \ldot \vk X(t_n))$ with respect to $Y$. It is well-known that exponential tilting of multivariate Gaussian df  is again a multivariate Gaussian, see for instance \cite[Lem 6.1]{Htilt}. The only thing that changes under the exponential tilting is the trend which can be calculated as in the aforementioned lemma for each component separately.\\
\QED

{ \bf Acknowledgement}:  I thank Krzys D\c{e}bicki,  Dmitry Korshunov, Ilya Molchanov and  Philippe Soulier for  discussions on related topics. A special great thanks goes to both referees for their very deep reviews, numerous suggestions and corrections as well as new ideas for continuations. Partial 	support by SNSF Grant 200021-196888 is kindly acknowledged.

%

\bibliographystyle{ieeetr}
\bibliography{EEEA}
\end{document}